\documentclass[a4paper,12pt]{article}
\usepackage{longtable}
\usepackage{setspace}
\usepackage{multirow}
\usepackage{epsfig}
\usepackage{amsmath}
\usepackage{amssymb}
\usepackage{tabularx}
\usepackage{enumerate}
\usepackage{graphicx}
\usepackage{subfigure}
\usepackage{cases}
\usepackage[compress]{cite}
\usepackage{lineno}
\usepackage{color}
 \usepackage{epstopdf}
\usepackage{geometry}
\newtheorem{thm}{Theorem}[section]
\newtheorem{lem}[thm]{Lemma}

\newtheorem{prop}[thm]{Proposition}
\newtheorem{rmk}[thm]{Remark}

\newtheorem{pppp}{Proof}

\newcommand{\qed}{\hspace{1em}\mbox{\raisebox{0.65ex}{\fbox{}}}}

\numberwithin{equation}{section}

\newcommand{\be}{\begin{equation}}
\newcommand{\ee}{\end{equation}}
\newcommand\bes{\begin{eqnarray}} \newcommand\ees{\end{eqnarray}}
\newcommand{\bess}{\begin{eqnarray*}}
\newcommand{\eess}{\end{eqnarray*}}
\newcommand{\bpf}{{\bf Proof.\ \ }}
\newcommand{\epf}{\mbox{}\hfill $\Box$}

\begin{document}

\thispagestyle{empty}
\title{Threshold dynamics of a nonlocal dispersal SIS epidemic model with free boundaries
\thanks{The first author is supported by the Postgraduate Research \& Practice
Innovation Program of Jiangsu Province (KYCX21-3188), the second
author is supported by the National Research Foundation of Korea(NRF) grant funded by the Korea government (MSIT) (NRF-2022R1F1A1063068)
and the third author is supported by the National Natural Science Foundation of China (Grant No. 12271470).}
}
\date{\empty}
\author{Yachun Tong$^{a}$, Inkyung Ahn$^{b}$ and Zhigui Lin$^a\thanks{Corresponding author.
Email: zglin@yzu.edu.cn (Z. Lin).}$\\
{\small $^a$ School of Mathematical Science, Yangzhou University, Yangzhou 225002, China}\\
{\small $^b$ Department of Mathematics, Korea University, Sejong 339-700, Republic of Korea}
}
\maketitle
\begin{quote}
\noindent
{\bf Abstract.} {
 \small
To study the influence of the moving front of the infected interval and the spatial
movement of individuals on the spreading or vanishing of infectious disease, we consider a
nonlocal SIS (susceptible-infected-susceptible) reaction-diffusion model with media coverage,
hospital bed numbers and free boundaries. The principal eigenvalue of the integral operator
is defined, and the impacts of the diffusion rate of infected individuals and
interval length on the principal eigenvalue are analyzed. Furthermore, sufficient conditions for
spreading and vanishing of the disease are derived.
Our results show that large media coverage and hospital bed numbers are beneficial
to the prevention and control of disease.  The difference between the model with nonlocal diffusion and that with local diffusion
 is also discussed and nonlocal diffusion leads to more possibilities.
}

\noindent {\it MSC:} 35K57, 92D30; secondary: 35R35.

\medskip
\noindent {\it Keywords:} SIS model; Free boundary; Nonlocal diffusion; Spreading and vanishing
\end{quote}

\section{Introduction}

With the emergence and outbreak of COVID-19 \cite{BTSK,M1} in recent years, infectious disease
models have become one of the most popular research topics. To study the spread and
dynamics of COVID-19, most scholars use the SIR (susceptible-infected-recovered) \cite{BTSK,LWLB},
SEIR (susceptible-exposed-infected-recovered) \cite{M1,LD} and SEAIR (susceptible-exposed-asymptomatic-infectious-removed) \cite{B,ZY} models to describe the spread of COVID-19.
Meanwhile, the classical SIS model has received great attention in mathematical epidemiology.

Considering the impact of the spatial heterogeneity of the environment and the movement of individuals on
infectious diseases, Allen et al. in \cite{AL} proposed and discussed an SIS reaction-diffusion
system
\begin{equation}
\left\{\begin{array}{lll}
S_{t}-d_S\Delta S=-\frac{\beta (x) SI}{S+I}+\gamma(x)I,\; &\, t>0,\ x\in\Omega, \\[2mm]
I_{t}-d_I\Delta I=\frac{\beta (x) SI}{S+I}-\gamma(x)I,\; &\, t>0,\ x\in\Omega, \\[2mm]
\frac{\partial S}{\partial\eta}=\frac{\partial I}{\partial\eta}=0,\; &\, t>0,\
x\in\partial\Omega.
\end{array} \right.
\label{aa}
\end{equation}
Here, $\Omega\subset\mathbb{R}^{n}$ $(n\geq1)$ is a bounded domain; $S(t,x)$ and $I(t,x)$
indicate the density of susceptible and infected individuals at location $x$ and time $t$,
respectively; $d_S$ and $d_I$ are positive constants that account for the diffusion rate
of susceptible and infected individuals, respectively; and the positive bounded H$\ddot{o}$lder
continuous functions $\beta(x)$ and $\gamma(x)$ can be interpreted as rates of disease
transmission and recovery for $x\in\Omega$, respectively. The authors in \cite{AL}
mainly discussed the existence, uniqueness and stability of DFE (disease-free equilibrium)
and EE (endemic equilibrium) and used the basic reproduction number ${\mathcal{R}}_0$ to
characterize the risk of the region. Afterwards, Peng and Liu \cite{PL} confirmed the conjecture
proposed by Allen et al. in \cite{AL} that a unique EE is globally asymptotically stable in
some special cases. Further results that the effect of individual movement (large or small)
on the existence and disappearance of disease were obtained in \cite{P}. For more results
of the SIS reaction-diffusion model, one can see \cite{HHL,PZ,WZ} and the references therein.

It is easy to find that the above articles are devoted to the study of SIS
models on a fixed domain. In real life, the movement of species leads to changes in
biological habitats, and in mathematics, the free boundary can be used to describe this
phenomenon, such as the healing of wounds \cite{CF} and the expansion of new species or invasive
species \cite{DL,CLW,LL,TN}. Free boundary problems can also be used to describe the
transmission of disease, such as the SIRS model \cite{CLWY}, SIS model \cite{GKLZ}, SIR model
\cite{HW,ZGL} and references therein.

To explore the moving front of the infected individual, Wang and Guo \cite{WG} introduced
the free boundary and studied the dynamics of the following SIS reaction-diffusion model:
\begin{equation}
\left\{ \begin{array}{llllll}
S_{t}-d\Delta S=\sigma-\mu S-\beta(x)SI+\gamma(x)I, &t>0,\, x\in \mathbb{R},\\[2mm]
I_{t}-d\Delta I=\beta(x)SI-\mu I-\gamma(x)I,& t>0,\, x\in(g(t),h(t)),\\[2mm]
I(t,x)=0, &t>0,\, x\in \mathbb{R}\backslash(g(t),h(t)),\\[2mm]
g'(t)=-k I_{x}(t,g(t)),\ g(0)=-h_{0},&t \geq 0, \\[2mm]
h'(t)=-k I_{x}(t,h(t)),\ h(0)=h_{0},&t\geq 0, \\[2mm]
S(0,x)=S_0(x),\ I(0,x)= I_0(x),&x\in \mathbb{R}.
\end{array}\right.
\label{a}
\end{equation}
The basic reproduction number was given, and the spreading-vanishing dichotomy was established.
Some conditions for disease spreading or vanishing were presented by investigating the
effect of the diffusion rate $(d)$, initial value ($I_0$) and expanding capability ($k$)
on the asymptotic behavior of the infected individuals.

It is widely known that random dispersal or local diffusion describes the local behavior
of the movements of organisms between adjacent spatial locations \cite{LSW}. Briefly, the
classical Laplace diffusion operator is used to describe that the movement of the infectious
agent and infected population only occurs between adjacent spatial positions \cite{XLR}.
However, Murray \cite{M} noted that a local or short-range diffusive flux proportional
to the gradient is not suitable to characterize some biological phenomena. In the real
world, the movements and interactions of some organisms occur at nonadjacent spatial
positions, and such dispersal is called nonlocal diffusion \cite{DWW}.
%%%
Nonlocal diffusion can occur naturally through dispersal and migration or facilitated by human activities. It can positively impact a population's genetic diversity and long-term viability, but it can also introduce disease or invasive species into new territories.
 Local diffusion usually involves individuals migrating short distances within a defined area, such as a habitat patch or a specific population. Various factors, including random movement, resource competition, and response to environmental conditions, can cause it.
On the other hand, nonlocal diffusion refers to the movement of individuals between different regions or subpopulations. This migration typically involves individuals moving long distances or crossing geographic barriers such as rivers or mountain ranges.
%%%

Recently, nonlocal diffusion equations have attracted extensive attention and have
been used to characterize long-range dispersal in population ecology \cite{CLW,LSW}.
In addition, scholars have extensively investigated infectious disease models
with nonlocal diffusion, such as the West Nile virus model \cite{DN}, SIS epidemic
model \cite{1,YLR}, and SIR reaction-diffusion model \cite{ZLC,FLY}. For other epidemic
models with nonlocal diffusion, see references \cite{WWW,WD1,CD} and references therein.

In addition, there are many factors that affect the spread of infectious disease, such as
the contact transmission rate and the recovery rate. Educating the public about the disease
through mass media (such as television, radio, newspapers, billboards, internet, magazines, etc.),
is one of the important precautions. Therefore, media coverage can indirectly reduce the
contact rate between people and infectious diseases, thus reducing the contact transmission
rate of infectious diseases \cite{RPZ}. In general, the main factor impacting the recovery
rate is the availability of health care (such as the number of physicians, nurses, hospital beds
and isolation places). In fact, health and medical institutions use the hospital bed-population
ratio (HBPR) (the number of hospital beds per 10000 people) as a method of reckoning available
resources to the public \cite{MM}.

Taking into account nonlocal diffusion, media coverage and hospital bed numbers, we consider
the following  nonlocal dispersal SIS epidemic model with a free boundary:
{\small \be\left\{
\begin{array}{ll}
S_{t}=d\mathcal{L}_1[S]+\sigma-\mu_{1}S-\beta(m(x),I,x) SI+\gamma(b(x),I,x)I,&t>0,\,
x\in \mathbb{R}, \\[2mm]
I_{t}=d\mathcal{L}_2[I;g,h]-\mu_{2}I+\beta(m(x),I,x)SI-\gamma(b(x),I,x)I, &t>0,\,
x\in(g(t), h(t)),\\[2mm]
I(t,x)=0,  &t\geq0,\ x\in\mathbb{R}\backslash(g(t),h(t)),\\[2mm]
h'(t)=k\int_{g(t)}^{h(t)}\int_{h(t)}^{+\infty}J(x-y)I(t,x)dydx, &t>0,\\[2mm]
g'(t)=-k\int_{g(t)}^{h(t)}\int_{-\infty}^{g(t)}J(x-y)I(t,x)dydx, &t>0,\\[2mm]
S(0,x)=S_{0}(x),\ g(0)=-h_{0},\ h(0)=h_{0},   &x\in\mathbb{R},\\[2mm]
I(0,x)= I_0(x),   &x\in(-h_0,h_0),
\end{array} \right.
\label{a02} \ee }
where
$$
\mathcal{L}_1[S]=\int_\mathbb{R}J(x-y)S(t,y)dy-S(t,x),
$$
$$
\mathcal{L}_2[I;g,h]=\int_{g(t)}^{h(t)}J(x-y)I(t,y)dy-I(t,x),
$$
and $d,\,S(t,x)$ and $I(t,x)$ have the same epidemiological interpretation as in \eqref{aa}.
The constants $\sigma,\,\mu_1$ and $\mu_2$ are positive, where $\sigma$ accounts for
the environment carrying capability; the natural mortality rate of the susceptible individuals
is expressed by $\mu_1$, and $\mu_2$ denotes the sum of the natural mortality and
disease-caused death rates of the infected individuals. The functions
$\beta(m(x),I,x),\,\gamma(b(x),I,x),\,m(x),\,b(x)$ are nonnegative, where $m(x)$
represents the media coverage, and $b(x)$ stands for the number of hospital beds.
In this paper, we assume that

$(1)$ the contact infectious rate $\beta(m(x),I,x)$ is Lipschitz continuous and monotonically
decreasing in $m(x)$ and increasing in $I$;

$(2)$ the recovery rate $\gamma(b(x),I,x)$ is Lipschitz continuous and increasing in $b(x)$ and
monotonically decreasing in $I$;

$(3)$ $\beta_{I}(m(x),I,x)$ and $\gamma_{I}(b(x),I,x)$ are continuous and bounded for $m(x)\in [0,\infty),I\in [0,\infty)$ and $x\in (-\infty, \infty)$.\\

For instance, Cui and Zhu \cite{CTZ} used the function $\beta(I)=\beta e^{mI}$ to
model the impact of media coverage on the transmission rate, and Shan and Zhu
\cite{SZ} used the function $\gamma(b,I,x)=\gamma_{0}+(\gamma_1-\gamma_0)\frac{b}{b+I}$
to describe the hospital resource impact factors.

Recalling that $S(t,x)$ denotes the density at point $x$ and time $t$, the kernel function
$J(x-y)$ is regarded as the probability distribution of jumping from place $y$ to place $x$,
then the integral operator $\int_\mathbb{R}J(x-y)S(t,y)dy$ accounts for the rate at which the
individuals are gathering at point $x$ from all other places, and $-S(t,x)$ is the rate at
which the individuals are leaving at point $x$ to other places. In addition, the infected
individuals stay in the infected interval $(g(t),h(t))$. We further suppose that the initial
function $S_{0}(x)$ satisfies
\be
S_{0}(x)\in C(\mathbb{R})\cap L^{\infty}(\mathbb{R})\, \ \textrm{and}\, \ S_{0}(x)>0\, \
\textrm{in}\, \,  \mathbb{R},
\label{a03}\ee
and $I_{0}(x)$ satisfies
\be
I_{0}(x)\in C([-h_0,h_0]),\, I_{0}(\pm h_0)=0,\, I_{0}(x)>0\, \ \textrm{in}\, \, (-h_0,h_0).
\label{a04}\ee

For system \eqref{a02}, we assume that the kernel function $J:\mathbb{R}\rightarrow\mathbb{R}$ is continuous and nonnegative and has the properties
$$
\mathbf{(J)}: J\in C(\mathbb{R})\cap L^{\infty}(\mathbb{R}) \, {\rm \, is \, symmetric},
\ J(0)>0,\ \int_{\mathbb{R}}J(x)dx=1.
$$
The free boundary conditions $h'(t)=k\int_{g(t)}^{h(t)}\int_{h(t)}^{+\infty}J(x-y)I(t,x)dydx$
and $g'(t)=-k\int_{g(t)}^{h(t)}\int_{-\infty}^{g(t)}J(x-y)I(t,x)dydx$ in \eqref{a02} imply
that the expanding rate of the interval $(g(t),h(t))$ is determined by the infected individuals
and is proportional to the outward flux of the infected individuals across the interval $(g(t),h(t))$ \cite{CDLL}.

It is worth mentioning that there are links and differences between local diffusion and
nonlocal diffusion. Local diffusion, expressed by the Laplace operator $\Delta u$ (the Laplace
in $\mathbb{R}^n$, $n\geq2$) or $u_{xx}$ (in one-dimensional space), is used to describe the
influence between adjacent positions, and nonlocal diffusion, expressed by the integral operator
(is given by $\int_\mathbb{R}J(x-y)u(t,y)dy-u(t,x)$), is used to describe long-distance
dispersal. However, the Laplace operator can be regarded as a local approximation of a nonlocal
diffusion operator. In fact, when $J(\cdot)$ is symmetric and has compact supports, such
as $J(x)=(1/\epsilon)K(x/\epsilon)$ with $0<\epsilon\ll1$ and $K(x)$ is a general mollification
function with support $x\in[-1,1]$, we can transform nonlocal operators into local operators
by using the Taylor formula \cite{LXZ}.

This article is organized as follows: the existence and uniqueness of the global solution are
given in Section 2. Section 3 is devoted to defining and studying the properties of the
principal eigenvalue. Section 4 gives some sufficient conditions for the disease to spread
or vanish. Finally, a brief discussion is presented in Section 5.

\section{Global existence and uniqueness}
In this section, we assume that $h_{0}>0$, $S_{0}(x)$ and $I_{0}(x)$ satisfy \eqref{a03}
and \eqref{a04}. For any given $T>0$, we first introduce the notations as follows:
\bess
\begin{array}{lll}
&\mathbb{H}_{T}:=\{h\in C([0, T]): h(0)=h_0,\, \ \inf\limits_{0\leq t_{1}<t_{2}\leq T}\frac{h(t_{2})-h(t_{1})}{t_{2}-t_{1}}>0\}, \\[2mm]
&\mathbb{G}_{T}:=\{g\in C([0, T]): -g\in \mathbb{H}_{T}\}, \\[2mm]
&D^{g,h}_{T}:=\{(t,x)\in\mathbb{R}^2: 0<t\leq T,\, \ g(t)<x<h(t)\}, \\[2mm]
&D^{h_{0}}_{T}:=\{(t,x)\in\mathbb{R}^2: 0<t\leq T,\, \ -h_{0}<x<h_{0}\}, \\[2mm]
&D^{\infty}_{T}:=\{(t,x)\in\mathbb{R}^2: 0<t\leq T,\, \ x\in\mathbb{R}\}, \\[2mm]
&X^{S_{0}}_{T}:=\{\phi(t,x)\in C(D^{\infty}_{T})\cap L^{\infty}(D^{\infty}_{T}):
\phi(0,x)=S_{0}(x)\, \ {\rm in}\, \ \mathbb{R},\, \ \phi(t,x)\geq0 \, \ {\rm in}\,
\ D^{\infty}_{T}\}, \\[2mm]
&X^{I_{0}}_{T}:=\{\psi(t,x)\in C(D^{\infty}_{T}): \psi(0,x)=I_{0}(x)\, \ {\rm in}\,
\ [-h_{0},h_{0}],\, \ \psi(t,x)\geq0 \, \ {\rm in}\, \ D^{g,h}_{T}, \\[2mm]
&\qquad \qquad \qquad \qquad \qquad \qquad \psi(t,x)=0\, \ {\rm for}\, \ t\in(0,T),\,
\ x\in \mathbb{R}\backslash(g(t),h(t))\}.
\end{array}
\eess

To prove the existence and uniqueness of the global solution of problem \eqref{a02},
we first give the following result for problem \eqref{a02} without a free boundary.

\begin{lem}
\label{lem2.1}
For any given $T>0$ and $(g,h)\in\mathbb{H}_{T}\times\mathbb{G}_{T}$, the problem
\be\left\{
\begin{array}{ll}
S_{t}=d\mathcal{L}_1[S]+\sigma-\mu_{1}S-\beta(m(x),I,x)SI+\gamma(b(x),I,x)I, &0<t\leq T,
\ x\in \mathbb{R}, \\[2mm]
I_{t}=d\mathcal{L}_2[I;g,h]-\mu_{2}I+\beta(m(x),I,x)SI-\gamma(b(x),I,x)I,  &0<t\leq T,
\ x\in(g(t), h(t)),\\[2mm]
I(t,x)=0,  &0\leq t\leq T,\ x\in\mathbb{R}\backslash(g(t),h(t)),\\[2mm]
S(0,x)=S_{0}(x), &x\in \mathbb{R}, \\[2mm]
I(0,x)= I_0(x),   &x\in(-h_0,h_0)
\end{array} \right.
\label{le2.1-1} \ee
admits a unique solution $(S_{g,h},I_{g,h})\in C(\overline{D}_{T}^{\infty})\times C(\overline{D}_{T}^{g,h})$. Moreover,
\be
0<S_{g,h}(t,x)\leq A\qquad  {\rm for}\, \ {\rm any} \,\,(t,x)\in D_{T}^{\infty},
\label{le2.1-2}
\ee
\be
0<I_{g,h}(t,x)\leq A\qquad  {\rm for}\, \ {\rm any} \,\,(t,x)\in D_{T}^{g,h},
\label{le2.1-3}
\ee
where $A=\max\{\frac{\sigma}{\mu_{1}}, \|S_{0}\|_{\infty}+\|I_{0}\|_{\infty}\}$.
\end{lem}

\bpf
The main idea of this proof comes from \cite{ZLC}. We divide the proof into three steps.

\textbf{Step 1.} The parameterized ODE problem.

For any given $x\in\mathbb{R}$, $s\in(0,T]$, denote
\bess
t_{x}=\left\{
\begin{array}{ll}
t_{x}^{g},\  &x\in(g(s),-h_0)\, \, \mbox{and}\, \, x=g(t_{x}^{g}),\\[2mm]
0,\  &x\in[-h_0,h_0],\\[2mm]
t_{x}^{h},\  &x\in(h_0,h(s))\, \, \mbox{and}\, \, x=h(t_{x}^{h}),\\[2mm]
s,\  &x\in\mathbb{R}\backslash(g(s),h(s)).
\end{array} \right.
\eess
Clearly, $t_{x}>0$ for $x\in\mathbb{R}\backslash[-h_0,h_0]$, $t_{x}<s$ for $x\in(g(s),h(s))$.
For any given
$(\phi,\psi)\in X_{s}^{S_0}\times X_{s}^{I_0}$, define
$$
A_{1}=\max\{A,\,\frac{\sigma+\|\psi\|_{\infty}\sup\gamma}{\mu_1},\,\|\phi\|_{\infty}\},\, \,
\, \, A_2=\max\{A,\,\frac{(d+A_1\sup\beta)\|\psi\|_{\infty}}{d+\mu_{2}}\}.
$$
We discuss it in the following two cases:

Case 1: $x\in\mathbb{R}\backslash[-h_0,h_0]$, $t\in[0,t_{x}]$.

Clearly, $I(t,x)=0$ for $(t,x)\in[0,t_{x}]\times\mathbb{R}\backslash[-h_0,h_0]$.
Consider the ODE problem
\begin{equation}
\left\{\begin{array}{lll}
S_{t}=d\int_\mathbb{R}J(x-y)\phi(t,y)dy-dS+\sigma-\mu_{1}S,\,\;&\,0<t\leq t_{x}, \\[2mm]
S(0,x)=S_{0}(x),\;&\, x\in\mathbb{R}\backslash[-h_0,h_0].
\end{array} \right.
\label{le2.1-4}
\end{equation}
For any $S_1,\,S_2\in[0,A_1]$,
$$
\begin{array}{llllll}
&&|d\int_\mathbb{R}J(x-y)\phi(t,y)dy-dS_1+\sigma-\mu_{1}S_1-d\int_\mathbb{R}J(x-y)
\phi(t,y)dy+dS_2
-\sigma+\mu_{1}S_2|\\[2mm]
&=&(d+\mu_{1})|S_1-S_2|.
\end{array}
$$
Therefore, $F:=d\int_\mathbb{R}J(x-y)\phi(t,y)dy-dS+\sigma-\mu_{1}S$ is Lipschitz continuous
in $S$ for $S\in[0,A_1]$. By the fundamental theory of ODEs, problem \eqref{le2.1-4} has a
unique solution $S_{\phi}(t,x)$ defined in $t\in[0,\widehat{t}_{x})$, and
$S_{\phi}(t,x)$ is continuous in both $t$ and $x$. To see that $t\rightarrow S(\cdot,x)$
can be uniquely extended to $[0,t_{x}]$, we need to prove that if $S_{\phi}(t,x)$ is uniquely
defined for $t\in[0,\widehat{t}_{x}]$ with $\widehat{t}_{x}\in(0,t_{x}]$, then
$$
0\leq S_{\phi}(t,x)\leq A_1,\, \, \, \mbox{for}\,\, \, t\in[0,\widehat{t}_{x}]\, \,
\mbox{and}\,\, x\in \mathbb{R}\backslash[-h_0,h_0].
$$
Obviously,
$$
\begin{array}{llllll}
&&d\int_{\mathbb {R}}J(x-y)\phi(t,y)dy-dA_{1}+\sigma-\mu_{1}A_{1}\\[2mm]
&\leq&d\|\phi\|_{\infty}-dA_{1}+\sigma-\mu_{1}A_{1}\\[2mm]
&\leq&0,
\end{array}
$$
and $\|S_{0}\|_{\infty}\leq A_{1}$. Thanks to the direct comparison argument, one can
derive $S_{\phi}(t,x)\leq A_1$ for $t\in[0,\widehat{{t}}_{x}]$ and $
x\in \mathbb{R}\backslash[-h_0,h_0]$. We use similar method to prove that
$S_{\phi}(t,x)\geq 0$ for $t\in[0,\widehat{t}_{x}],\, x\in \mathbb{R}\backslash[-h_0,h_0]$.

Case 2: $x\in(g(s),h(s))$, $t\in[t_{x},s]$.

Define
$$
\widehat{S}_{\phi}(x)=\left\{
\begin{array}{ll}
S_{0}(x),\  &x\in[-h_0,h_0]\\[2mm]
S_{\phi}(t_{x},x),\  &x\notin[-h_0,h_0]
\end{array} \right.
\,\,\, \mbox{and}\,\,\,\
\widehat{I}(x)=\left\{
\begin{array}{ll}
I_{0}(x),\  &x\in[-h_0,h_0]\\[2mm]
0,\  &x\notin[-h_0,h_0].
\end{array} \right.
$$
Consider the ODE problem
\begin{equation}
\left\{\begin{array}{lll}
S_{t}=F_{1}(t,x,S,I),\,\;&\,t_{x}<t\leq s, \\[2mm]
I_{t}=F_{2}(t,x,S,I),\;&\, t_{x}<t\leq s,\\[2mm]
S(t_{x},x)=\widehat{S}_{\phi}(x),\, \, I(t_{x},x)=\widehat{I}(x),\;&\,x\in(g(s),h(s))
\end{array} \right.
\label{le2.1-5}
\end{equation}
with
$$
F_{1}=d\int_\mathbb{R}J(x-y)\phi(t,y)dy-dS+\sigma-\mu_{1}S+\gamma(b,I,x)\psi-\beta(m,I,x) SI,
$$
$$
F_{2}=d\int_{g(t)}^{h(t)}J(x-y)\psi(t,y)dy-dI-\mu_{2}I-\gamma(b,I,x)I+\beta(m,I,x) S\psi.
$$
For any $(S_{i},I_{i})\in[0,A_{1}]\times[0,A_{2}]$($i=1,\,2$), obviously,
$F_{i}(t,x,S,I)$ is Lipschitz continuous in $(S,I)$ for $(S_{i},I_{i})\in[0,A_{1}]\times[0,A_{2}]$
by the continuity and monotonicity of $\beta(m(x),I,x)$ and $\gamma(b(x),I,x)$, and it
is uniformly continuous for $x\in(g(s),h(s))$ and $t\in[t_{x},s]$. In addition, $F_{i}(t,x,S,I)$
is continuous in all its variables in this range. Problem \eqref{le2.1-5} has a unique
solution $(S_{\phi,\psi}(t,x),\,I_{\phi,\psi}(t,x))$ for $t\in[t_{x},s_{x})$, and $(S_{\phi,\psi}(t,x),\,I_{\phi,\psi}(t,x))$ is continuous in both $t$ and $x$ by the
fundamental theorem of ODEs.

To show that $(S_{\phi,\psi}(t,x),\,I_{\phi,\psi}(t,x))$ can be uniquely extended to $[t_x,s]$,
it suffices to prove that if $(S_{\phi,\psi}(t,x),\,I_{\phi,\psi}(t,x))$ is uniquely
defined for $t\in[t_{x},\widehat{t}]$ with $\widehat{t}\in(t_{x},s]$, then
\be
0\leq S_{\phi,\psi}(t,x)\leq A_{1},\,\, 0\leq I_{\phi,\psi}(t,x)\leq A_{2}\, \,\,\, \mbox{for}\,\,\,t\in[t_{x},\widehat{t}].
\label{le2.1-5.1}
\ee
In fact, it is easy to see that
$$
\begin{array}{llllll}
&&F_{1}(t,x,A_{1},A_{2})\\[2mm]
&=&d\int_\mathbb{R}J(x-y)\phi(t,y)dy-dA_{1}+\sigma-\mu_{1}A_{1}
+\gamma(b,A_{2},x)\psi-\beta(m,A_2,x) A_{1} A_{2}\\[2mm]
&\leq&d\|\phi\|_{\infty}-dA_{1}+\sigma-\mu_{1}A_{1}
+\gamma(b,A_{2},x)\|\psi\|-\beta(m,A_2,x) A_{1}A_2\\[2mm]
&<&d\|\phi\|_{\infty}-dA_{1}+\sigma-\mu_{1}A_{1}
+\|\psi\|_{\infty}\sup\gamma\\[2mm]
&\leq&0
\end{array}
$$
and
$$
\begin{array}{llllll}
&&F_{2}(t,x,A_{1},A_{2})\\[2mm]
&=&d\int_{g(t)}^{h(t)}J(x-y)\psi(t,y)dy-dA_{2}-\mu_{2}A_{2}
-\gamma(b,A_{2},x)A_{2}+\beta(m,A_2,x) A_{1}\psi\\[2mm]
&\leq&d\int_{g(t)}^{h(t)}J(x-y)\psi(t,y)dy-dA_{2}-\mu_{2}A_{2}+A_{1}\|\psi\|_{\infty}
\sup\beta\\[2mm]
&\leq&0.
\end{array}
$$
Since $A_{1}\geq\|S_{0}\|_{\infty},\,A_{2}\geq\|I_{0}\|_{\infty}$, we have
$S_{\phi,\psi}(t,x)\leq A_{1}$ and $I_{\phi,\psi}(t,x)\leq A_{2}$ in $t\in[t_{x},\widehat{t}]$
by the comparison argument. The left part of \eqref{le2.1-5.1} can be obtained
similarly by using $F_{i}(t,x,0,0)\geq0$ ($i=1,\,2$).

\textbf{Step 2.} A fixed point theorem.

For any $s\in(0,T)$, we note
$$
X_{s}^{S_{0}}:=\{\phi|_{\overline{D}_{s}^{\infty}}:\phi\in X_{T}^{S_{0}}\},\,\,\,
X_{s}^{I_{0}}:=\{\psi|_{\overline{D}_{s}^{g,h}}:\psi\in X_{T}^{I_{0}}\}.
$$
Denote
\bess
(\widehat{S}(t,x),\widehat{I}(t,x))=\left\{
\begin{array}{ll}
(S_{\phi(t,x)},0),\  &x\in\mathbb{R}\backslash[-h_{0},h_{0}],\, \, t=[0,t_{x}],\\[2mm]
(S_{\phi,\psi}(t,x),I_{\phi,\psi}(t,x)),\  &x\in(g(s),h(s)),\, \, t\in[t_{x},s],
\end{array} \right.
\eess
where $S_{\phi}(t,x),\,S_{\phi,\psi}(t,x)$ and $I_{\phi,\psi}(t,x))$ are given in Step 1.
By Step 1, for any $(\phi,\,\psi)$, we have a unique solution $(\widehat{S},\widehat{I})$
for $t\in[0,s]$. It is easy to check that $\widehat{S}(t,x)$ is continuous in
$\overline{D}_{s}^{\infty}$, and $\widehat{I}(t,x)$ is continuous in $\overline{D}_{s}^{g,h}$
due to the continuous dependence of the ODE solution on the parameters. Therefore, $(\widehat{S},\widehat{I})\in X_{s}^{S_{0}}\times X_{s}^{I_{0}}$. Note that $X_{s}^{S_{0}}$ and $X_{s}^{I_{0}}$ are
complete metric spaces, respectively, with the norms
$$
d_{1}(\phi_1,\phi_2)=\|\phi_1-\phi_2\|_{C(\overline{D}_{s}^{\infty})},\, \, \, d_{2}(\psi_1,\psi_2)=\|\psi_1-\psi_2\|_{C(\overline{D}_{s}^{g,h})}.
$$
Hence, we find a mapping $\Gamma:X_{s}^{S_{0}}\times X_{s}^{I_{0}}\rightarrow X_{s}^{S_{0}}\times X_{s}^{I_{0}}$ by $\Gamma(\phi,\psi)=(\widehat{S},\widehat{I})$.

Setting
$$
M_1=\max\{A,\,4\|S_{0}\|_{\infty},\,\frac{4\sigma}{\mu_{1}},\,
\frac{4(\sigma+M_2)}{\mu_{1}+d}\},\,
M_2=\max\{A,\,2\|I_{0}\|_{\infty}\}.
$$
Define
$$
X_{s}^{M_1}=\{\phi|\,\phi\in X_{s}^{S_{0}},\,\,\|\phi\|_{C(\overline{D}_{s}^{\infty})}\leq M_1\},
$$
$$
X_{s}^{M_2}=\{\psi|\,\psi\in X_{s}^{I_{0}},\,\,\|\psi\|_{C(\overline{D}_{s}^{g,h})}\leq M_2\}.
$$
Using the same arguments as Lemma 2.1 in \cite{ZLC}, we can deduce that $\Gamma$ is a
contraction map and has a unique fixed point $(S^*,\,I^*)\in X_{s}^{M_1}\times X_{s}^{M_2}$
 for any $s\in(0,\widehat{s}]$ by the contraction mapping theorem, where $\widetilde{s}$
 relies on $d,\,M_1,\, \beta,\,\gamma$ and $M_2$.

To prove that $(S^*,\,I^*)$ is the unique solution \eqref{le2.1-1} for $t\in[0,s]$
with $s\in(0,\widehat{s}]$, it suffices to discuss that any nonnegative solution $(S,\,I)$
of \eqref{le2.1-1} for $t\in[0,s]$ belongs to $X_{s}^{M_1}\times X_{s}^{M_2}$.

We claim that
\be
S+I\leq A \, \,\,\, \mbox{for}\,\,\,t\in[0,s]\, \, \, \, \mbox{and}\, \, \, x\in\mathbb{R},
\label{le2.1-5.2}
\ee
which implies that
$$
0\leq S(t,x)\leq A, \,\,\,\,(t,x)\in[0,s]\times\mathbb{R},
$$
$$
0\leq I(t,x)\leq A, \,\,\,\,(t,x)\in[0,s]\times[g(t),h(t)].
$$
Consequently, we obtain that for any $s\in(0,\widehat{s}]$, \eqref{le2.1-1} admits a unique
solution for $t\in[0,s]$. To complete the proof, it only needs to prove that
claim \eqref{le2.1-5.2} is true.

Let $N=S+I$, then for $t\in[0,s]$ and $x\in(g(t), h(t))$,
$$\begin{array}{lll}
N_{t}&\leq&
d\int_{\mathbb{R}}J(x-y)N(t,y)dy-dN(t,x)-d(\int_{-\infty}^{g(t)}+\int_{h(t)}^{\infty})J(x-y)
I(t,y)dy+\sigma-\mu_{1}N\\[2mm]
&\leq&d\int_{\mathbb{R}}J(x-y)N(t,y)dy-dN(t,x)-\mu_{1}N+\sigma.
\end{array}$$
Since $\|S_{0}\|_{\infty}+\|I_{0}\|_{\infty}\leq A$, we have $N(t,x)\leq A$ for $t\in[0,s]$
and $x\in(g(t), h(t))$ by using the comparison principle.

While $t\in[0,s]$ and $x\in\mathbb{R}\backslash(g(t), h(t))$, then $I(t,x)=0$, which implies
$$
N_{t}=d\int_{\mathbb{R}}J(x-y)N(t,y)dy-dN(t,x)+\sigma-\mu_{1}N.
$$
It is clear that $N\leq A$ for $t\in[0,s]$ and $x\in\mathbb{R}\backslash(g(t), h(t))$.
Next, we prove that \eqref{le2.1-5.2} holds. We argue by contradiction and suppose that $\max\limits_{(t,x)\in[0,\,s]\times\mathbb{R}}N(t,x)>A$, there exists a point
$(t_0,x_0)\in[0,s]\times\mathbb{R}$ such that $\max N=N(t_0,x_0)>A$. According to the
above analysis, we can obtain that $x_0=g(t_0)$ or $x_0=h(t_0)$. Without loss of generality,
we assume that $x_0=g(t_0)$. Since $I(t_0,g(t_0))=0$, $S(t_0,x_0)$ satisfies
$$
S_{t}(t_0,x_0)=d\int_{\mathbb{R}}J(x_0-y)S(t_0,y)dy-dS(t_0,x_0)+\sigma-\mu_{1}S(t_0,x_0).
$$
Obviously, $S_t(t_0,x_0)\geq 0$ and $S(t_0,x_0)\leq A$, which contradicts the assumption:
$\max\limits_{(t,x)\in[0,\,s]\times\mathbb{R}}N(t,x)>A$.

\textbf{Step 3.} Extension of the solution.

We now prove that the unique solution of \eqref{le2.1-1} for $0<t\leq s$ can be extended
to $0<t\leq T$. In Step 2, $\widehat{s}$ depends only on $d,\,\gamma,\,\beta$, and  $A$. With
the help of the iterative method, we obtain that problem \eqref{le2.1-1} has a unique
solution for $t\in[0,T]$. We omit it here; see Step 3 of the proof in Lemma
2.1 in \cite{ZLC} for more details.
\epf

\begin{thm}
\label{thm2.2}
Assume that $\mathbf{(J)}$ holds. For any $S_{0}$ satisfying \eqref{a03} and $I_{0}$
satisfying \eqref{a04}, problem \eqref{a02} admits a unique positive solution
$(S(t,x),\, I(t,x);\, g(t),\, h(t))$ defined for all $t>0$.
\end{thm}
\bpf
We will prove this result by using Lemma \ref{lem2.1} and the fixed point theorem. For
any given $T>0$ and $(g^*,h^*)\in\mathbb{H}_{T}\times\mathbb{G}_{T}$, we know that
\eqref{le2.1-1} with $(g,h)=(g^*,h^*)$ has a unique solution $(S^*,I^*)$. Define
\begin{equation}
\left\{\begin{array}{ll}
\widetilde{g}=-h_0-k\int_{0}^{t}\int_{g^*(\tau)}^{h^*(\tau)}\int_{-\infty}^{g^*(\tau)}
J(x-y)I^*(\tau,x)dydxd\tau, \\[2mm]
\widetilde{h}=h_0+k\int_{0}^{t}\int_{g^*(\tau)}^{h^*(\tau)}\int_{h^*(\tau)}^{+\infty}
J(x-y)I^*(\tau,x)dydxd\tau.
\end{array} \right.
\label{thm2.2-1}
\end{equation}
In view of $\mathbf{(J)}$ and $J(0)>0$, there exist constants $\epsilon_0\in(0,h_0/4)$
and $\delta_0>0$ such that
$$
J(x)\geq\delta_0\, \, \, \mbox{if} \, \, \, |x|\leq\epsilon_0.
$$
By virtue of the above inequality and proof of Theorem 2.1 in \cite{CDLL}, there exists
$$T_0=T_0(k,A,h_0,\epsilon_0,I_0,J)>0,$$
such that, for any $T\in(0,T_0]$,
$$
\sup\limits_{0\leq t_1<t_2\leq T}\frac{\widetilde{g}(t_2)-\widetilde{g}(t_1)}{t_2-t_1}
\leq-k\eta_1,\,\,
\inf\limits_{0\leq t_1<t_2\leq T}\frac{\widetilde{h}(t_2)-\widetilde{h}(t_1)}{t_2-t_1}
\geq k\eta_2,
$$
$$
\widetilde{h}(t)-\widetilde{g}(t)\leq2h_0+\epsilon_0/4\, \, \, \mbox{for}\, \, t\in[0,T],
$$
where
$$
\eta_1=\frac{1}{4}\epsilon_0\delta_0e^{-(d+\mu_2+\sup\gamma)T_0}
\int_{-h_0}^{-h_0+\frac{\epsilon_0}{4}}I_0(x)dx,\, \, \,
\eta_2=\frac{1}{4}\epsilon_0\delta_0e^{-(d+\mu_2+\sup\gamma)T_0}
\int_{h_0-\frac{\epsilon_0}{4}}^{h_0}I_0(x)dx.
$$
Let
\bess
\begin{array}{lll}
\Sigma_T:=\{(g,h)\in\mathbb{H}_{T}^{h_0}\times\mathbb{G}_{T}^{h_0}:
\sup\limits_{0\leq t_1<t_2\leq T}\frac{g(t_2)-g(t_1)}{t_2-t_1}\leq-k\eta_1,\\[2mm]
\qquad \qquad \inf\limits_{0\leq t_1<t_2\leq T}\frac{h(t_2)-h(t_1)}{t_2-t_1}\geq k\eta_2,
\, h(t)-g(t)\leq2h_0+\frac{\epsilon_0}{4}\, \, \, {\rm for}\, \ t\in[0,T]\},
\end{array}
\eess
and define the mapping $\mathcal{F}(g^*,h^*)=(\widetilde{g},\widetilde{h})$.
Clearly, %It is easy to see that
 the above analysis implies that
$$
\mathcal{F}(\Sigma_T)\subset\Sigma_T\, \, \rm{for}\, \, T\in(0,T_0].
$$
In the following, similar to the proof of Theorem 1.1 in \cite{ZLC}, we first prove that
$\mathcal{F}$ is a contraction mapping, and then $\mathcal{F}$ admits a unique fixed
point in $\Sigma_T$ by the contraction mapping theorem. Next, we can derive that
$(g,h)\in\Sigma_T$ holds for any solution $(S,I;g,h)$ of \eqref{a02} for $t\in[0,T]$,
that is, $(S,I;g,h)$ is a unique solution of \eqref{a02} for $t\in[0,T]$. Finally, we
can show that the solution $(S,I;g,h)$ of \eqref{a02} is uniquely extended to
$t\in(0,+\infty)$. Here, we omit the proof here; see Step 3 in the proof of Theorem 1.1 in
\cite{ZLC} for more details.

\epf

\section{The eigenvalue problem}
For any $-\infty<L_1< L_2<+\infty$ and $d>0$, denote
$$
\mathcal{L}_{\{(L_1,\,L_2),\,d\}}[\phi](x)=d\int_{L_1}^{L_2}J(x-y)\phi(y)dy-d\phi(x).
$$
Now, we introduce some results on the principal eigenvalue of the linear operator
$\mathcal{L}_{\{(L_1,\,L_2),\,d\}}+a(x): C([L_1,L_2])\mapsto C([L_1,L_2])$
defined by
$$
(\mathcal{L}_{\{(L_1,\,L_2),\,d\}}+a(x))[\phi](x)=d\int_{L_1}^{L_2}J(x-y)\phi(y)dy
-d\phi(x)+a(x)\phi(x),
$$
where $a(x)=\frac{\sigma\beta(m(x),0,x)}{\mu_{1}}-\mu_{2}-\gamma(b(x),0,x)\in C([L_1, L_2])$
and $J$ satisfies $\mathbf{(J)}$. Furthermore, we assume
%$$
%\mathbf{(H)}: a(x)\, {\rm \, is \, Lipschitz \, continuous \, and \, achieves \, its \, maximum \, in} \, [L_1, L_2] \, {\rm \, at \, some \, point} \, x_0\in(L_1, L_2).
%$$

\vspace{0.2in}
$\mathbf{(H)}: a(x)$ is Lipschitz continuous and  achieves its maximum in $[L_1, L_2]$  at some point  $x_0 \in(L_1, L_2)$.
\vspace{0.2in}

Define the generalized principal eigenvalue as
\be
\begin{array}{ll}
&\lambda_{p}(\mathcal{L}_{\{(L_1,\,L_2),\,d\}}+a(x))\\
&:=\inf\{\lambda\in \mathbb{R}\,|\, \exists\phi\in C([L_1,L_2]),\, \phi>0\ {\rm\  s.\ t.} \
(\mathcal{L}_{\{(L_1,\,L_2),\,d\}}+a(x))[\phi]\leq \lambda\phi {\rm\ in} \ (L_1,L_2)\}.
\end{array}
\label{c01}
\ee
Furthermore, we call it a principal eigenvalue if $\lambda_{p}(\mathcal{L}_{\{(L_1,\,L_2),\,d\}}+a(x))$ is an eigenvalue of the operator $\mathcal{L}_{\{(L_1,\,L_2),\,d\}}+a(x)$ with a continuous and positive eigenfunction. Recalling that $a(x)$ is Lipschitz continuous and achieves a global maximum in $(L_1,L_2)$, $a(x)$ automatically satisfies the condition $\frac{1}{(\sup_{x\in(L_1,\,L_2)}a(x))-a(x)}\not\in L^{1}$. Therefore, it follows from Theorem 1.1 or Theorem 1.2 in \cite{C} that the generalized principal eigenvalue $\lambda_{p}(\mathcal{L}_{\{(L_1,\,L_2),\,d\}}+a(x))$ is a principal eigenvalue.

In this section, we are interested in the properties of the generalized principal eigenvalue $\lambda_{p}(\mathcal{L}_{\{(L_1,\,L_2),\,d\}}+a(x))$ with media coverage $m(x)$ and hospital bed number $b(x)$ and asymptotic behavior of the principal eigenvalue in large and small interval lengths $(L_1,L_2)$ or diffusion rate $d$.

Before stating our primary result, we recall a useful proposition from \cite{C}.

\begin{prop}({\cite{C}})
\label{prop3.1}
The following assertions hold:

$(i)$ Assume $(L_{1}, L_{2})\subset(L_3, L_4)$. Then,
$$
\lambda_{p}(\mathcal{L}_{\{(L_1,\,L_2),\,d\}}+a(x))\leq\lambda_{p}(\mathcal{L}_{\{(L_{3},\,L_{4}),\,d\}}+a(x)).
$$

$(ii)$ Fix $L_1, L_2$ and suppose that $a_{1}(x)\leq a_{2}(x)$.  Then,
$$
\lambda_{p}(\mathcal{L}_{\{(L_1,\,L_2),\,d\}}+a_{1}(x))\leq
\lambda_{p}(\mathcal{L}_{\{(L_1,\,L_2),\,d\}}+a_{2}(x)).
$$
Moreover, if $a_{1}(x)+\delta<a_{2}(x)$ for some $\delta>0$, then
$$
\lambda_{p}(\mathcal{L}_{\{(L_1,\,L_2),\,d\}}+a_{1}(x))<
\lambda_{p}(\mathcal{L}_{\{(L_1,\,L_2),\,d\}}+a_{2}(x)).
$$

$(iii)$ $\lambda_{p}(\mathcal{L}_{\{(L_1,\,L_2),\,d\}}+a(x))$ is Lipschitz
continuous in $a(x)$. More precisely,
$$ |\lambda_{p}(\mathcal{L}_{\{(L_1,\,L_2),\,d\}}+a_{1}(x))-
\lambda_{p}(\mathcal{L}_{\{(L_1,\,L_2),\,d\}}+a_{2}(x))|\leq \|a_{1}(x)-a_{2}(x)\|_{\infty}.
$$
\end{prop}

Let us now analyze the impact of media coverage $m(x)$ and hospital bed number $b(x)$ on
the generalized principal eigenvalue. Obviously, the following result holds by Proposition
\ref{prop3.1} $(ii)$.

\begin{thm}
\label{thm3.2-1}
Suppose that $\mathbf{(J)}$ and $\mathbf{(H)}$ hold. Then,

$(i)$ $\lambda_{p}(\mathcal{L}_{\{(L_1,\,L_2),\,d\}}+a(x))$ is strictly monotone
decreasing in $m(x)$.

$(ii)$ $\lambda_{p}(\mathcal{L}_{\{(L_1,\,L_2),\,d\}}+a(x))$ is strictly monotone
decreasing in $b(x)$.
\end{thm}

From now on, we discuss the effect of interval length on the principal eigenvalue $\lambda_{p}(\mathcal{L}_{\{(L_1,L_2),\,d\}}+a(x))$.

\begin{thm}
\label{thm3.2}
Suppose that $\mathbf{(J)}$ and $\mathbf{(H)}$ hold;
then, the following three conclusions hold:

$(i)$ $\lambda_{p}(\mathcal{L}_{\{(L_1,\,L_2),\,d\}}+a(x))$ is continuous for $L_1, L_2\in(-\infty,+\infty)$.

$(ii)$ $\lim\limits_{L_1,\,L_2\to 0}\lambda_{p}(\mathcal{L}_{\{(L_1,\,L_2),\,d\}}+a(x))=a(0)-d$.

$(iii)$ $\lim\limits_{-L_1, L_2\to+\infty}\lambda_{p}(\mathcal{L}_{\{(L_1,\,L_2),\,d\}}+a(x))=\sup\limits_{x\in\mathbb{R}}a(x)$.
\end{thm}
\bpf
The proof of $(i)$ is similar to Proposition 3.4 in \cite{CDLL}, and we omit it here.

$(ii)$ Due to the continuity of $a(x)$, for any given $\epsilon>0$, there exists $h>0$ small enough such that
$$
|a(x)-a(0)|<\epsilon, \, \, x\in[-h,h].
$$

Since $\lambda_{p}(\mathcal{L}_{\{(-h,\,h),\,d\}}+a(x))$ is a principal eigenvalue, there exists a positive function $\phi(x)\in C([-h,h])$ such that
$$
d\int_{-h}^{h}J(x-y)\phi(y)dy-d\phi(x)+a(x)\phi(x)=\lambda_{p}(\mathcal{L}_{\{(-h,\,h),\,d\}}+a(x))\phi(x),
\,\,x\in[-h,h],
$$
which gives by integrating,
$$
\begin{array}{llllll}
&&|\lambda_{p}(\mathcal{L}_{\{(-h,\,h),\,d\}}+a(x))-a(0)+d|\\[2mm]
&=&|\frac{d\int_{-h}^{h}\int_{-h}^{h}J(x-y)\phi(y)\phi(x)dydx}{\int_{-h}^{h}\phi^2(x)dx}
+\frac{\int_{-h}^{h}a(x)\phi^2(x)dx}{\int_{-h}^{h}\phi^2(x)dx}-a(0)|\\[2mm]
&=&|\frac{d\int_{-h}^{h}\int_{-h}^{h}J(x-y)\phi(y)\phi(x)dydx}{\int_{-h}^{h}\phi^2(x)dx}
+\frac{\int_{-h}^{h}(a(x)-a(0))\phi^2(x)dx}{\int_{-h}^{h}\phi^2(x)dx}|\\[2mm]
&\leq&\frac{d\|J\|_{\infty}(\int_{-h}^{h}\phi(x)dx)^2}{\int_{-h}^{h}\phi^2(x)dx}+
|\frac{\int_{-h}^{h}(a(x)-a(0))\phi^2(x)dx}{\int_{-h}^{h}\phi^2(x)dx}|\\[2mm]
&\leq&2d\|J\|_{\infty}h+\epsilon\rightarrow\epsilon\, \, \,  \mbox{as}\, \, \, h\rightarrow0^{+}.
\end{array}
$$
From the arbitrariness of $\epsilon$, we have
$$
|\lambda_{p}(\mathcal{L}_{\{(-h,\,h),\,d\}}+a(x))-a(0)+d|\rightarrow0 \, \, \,  \mbox{as}\, \, \, h\rightarrow0^{+},
$$
which together with the continuity of $\lambda_{p}(\mathcal{L}_{\{(L_1,\,L_2),\,d\}}+a(x))$ about $L_1,\,L_2$ give
$$
\lim\limits_{L_1,\,L_2\to 0}\lambda_{p}(\mathcal{L}_{\{(L_1,\,L_2),\,d\}}+a(x))=a(0)-d.
$$

$(iii)$ According to the monotonicity of $\lambda_{p}(\mathcal{L}_{\{(L_1,\,L_2),\,d\}}+a(x))$ with
respect to interval $(L_1,L_2)$ and function $a(x)$ yields
$$
\lambda_{p}(\mathcal{L}_{\{(L_1,\,L_2),\,d\}}+a(x))\leq\lambda_{p}(\mathcal{L}_{\{(L_1,\,L_2),\,d\}}+
\sup\limits_{x\in\mathbb{R}}a(x))\leq\lambda_{p}(\mathcal{L}_{\{(-\infty,\,+\infty),\,d\}}+
\sup\limits_{x\in\mathbb{R}}a(x)).
$$
Consider the following eigenvalue problem
\begin{equation}
d\int_{-\infty}^{\infty}J(x-y)\phi(y)dy-d\phi(x)+
\phi(x)\sup\limits_{x\in\mathbb{R}}a(x)=
\lambda_{p}(\mathcal{L}_{\{(-\infty,\,+\infty),\,d\}}+\sup\limits_{x\in\mathbb{R}}a(x))\phi,
\,\,x\in\mathbb{R}.
\label{th3.2-1}
\end{equation}
It is easily seen that $\lambda_{p}(\mathcal{L}_{\{(-\infty,\,+\infty),\,d\}}+\sup\limits_{x\in\mathbb{R}}a(x))
=\sup\limits_{x\in\mathbb{R}}a(x)$. So the principal eigenvalue $\lambda_{p}(\mathcal{L}_{\{(L_1,\,L_2),\,d\}}+a(x))\leq\sup\limits_{x\in\mathbb{R}} a(x)$.
Therefore,
$$
\limsup\limits_{-L_1,\,L_2\to+\infty}\lambda_{p}(\mathcal{L}_{\{(L_1,\,L_2),\,d\}}+a(x))
\leq\sup\limits_{x\in\mathbb{R}}a(x).
$$

To prove $(iii)$, it suffices to prove that $\liminf\limits_{-L_1,\, L_2\to+\infty}\lambda_{p}(\mathcal{L}_{\{(L_1,\,L_2),\,d\}}+a(x))\geq\sup\limits_{x\in\mathbb{R}}a(x)$ holds. In fact, by the continuity of $a(x)$ and the definition of sup, for given $\epsilon>0$, there exists some $x_0\in\mathbb{R}$ such that
$$
\sup\limits_{x\in\mathbb{R}}a(x)-\epsilon\leq a(x_0).
$$
By $\mathbf{(J)}$, for given $\epsilon>0$, there exist $L_1<x_0-1$ and $L_2>x_0+1$ such that
$$
\int_{L_1}^{L_2}J(z)dz>1-\epsilon.
$$

Now take
$$
\delta_{n}(x-x_0)=\left\{
\begin{array}{ll}
k_{1}e^{-1/(1-n^2(x-x_0)^2)}>0,\  &x\in(x_0-1/n,x_0+1/n),\\[2mm]
=0,\  &x\notin(x_0-1/n,x_0+1/n),
\end{array} \right.
$$
where $k_{1}$ is positive and satisfies $\int_\mathbb{R}\delta_{n}(x-x_0)dx=1$. It is easy to check that
the sequence $\{\delta_{n}(x-x_0)\}$ weakly converges to some $\delta(x-x_0)$ in $L^{1}((L_1,L_2))$.
By the definition of $\lambda_{p}(\mathcal{L}_{\{(L_1,\,L_2),\,d\}}+a(x))$, one can easily obtain
\begin{equation}
\begin{array}{llllll}
&&d\int_{L_1}^{L_2}J(x-y)\delta_{n}(y-x_0)dy-d\delta_{n}(x-x_0)+a(x)\delta_{n}(x-x_0) \\[2mm]
&\leq&\lambda_{p}(\mathcal{L}_{\{(L_1,\,L_2),\,d\}}+a(x))\delta_{n}(x-x_0),
\,\,x\in(L_1,L_2).
\label{th3.2-3}
\end{array}
\end{equation}
Integrating the equation of \eqref{th3.2-3} over $(L_1,L_2)$ yields
$$
\begin{array}{llllll}
\lambda_{p}(\mathcal{L}_{\{(L_1,\,L_2),\,d\}}+a(x))&\geq&\frac{d\int_{L_1}^{L_2}\int_{L_1}^{L_2}J(x-y)
\delta_{n}(y-x_0)dydx-d\int_{L_1}^{L_2}\delta_{n}(x-x_0)dx+\int_{L_1}^{L_2}a(x)\delta_{n}(x-x_0)dx}
{\int_{L_1}^{L_2}\delta_{n}(x-x_0)dx}\\[2mm]
&\geq&\frac{d\int_{L_1}^{L_2}\int_{x_{0}-1}^{x_{0}+1}J(x-y)
\delta_{n}(y-x_0)dydx-d\int_{L_1}^{L_2}\delta_{n}(x-x_0)dx+\int_{L_1}^{L_2}a(x)\delta_{n}(x-x_0)dx}
{\int_{L_1}^{L_2}\delta_{n}(x-x_0)dx}\\[2mm]
&\geq&\frac{d\int_{x_{0}-1}^{x_{0}+1}\delta_{n}(y-x_0)[\int_{L_{1}-x_{0}+1}^{L_{2}-x_{0}-1}J(z)
dz]dy-d\int_{L_1}^{L_2}\delta_{n}(x-x_0)dx+\int_{L_1}^{L_2}a(x)\delta_{n}(x-x_0)dx}
{\int_{L_1}^{L_2}\delta_{n}(x-x_0)dx}\\[2mm]
&\geq&\frac{(d(1-\epsilon)-d)\int_{L_1}^{L_2}\delta_{n}(x-x_0)dx+\int_{L_1}^{L_2}a(x)\delta_{n}(x-x_0)dx}
{\int_{L_1}^{L_2}\delta_{n}(x-x_0)dx}\\[2mm]
&=&-d\epsilon+\frac{\int_{L_1}^{L_2}a(x)\delta_{n}(x-x_0)dx}
{\int_{L_1}^{L_2}\delta_{n}(x-x_0)dx},
\end{array}
$$
where we have used that $\int_{x_0-1}^{x_0+1}\delta_{n}(x-x_0)dx=\int_{L_1}^{L_2}\delta_{n}(x-x_0)dx$ for sufficiently large $n$. Therefore, by taking $n\rightarrow+\infty$,
$$
\begin{array}{llllll}
\liminf\limits_{-L_1,\,L_2\to+\infty}\lambda_{p}(\mathcal{L}_{\{(L_1,L_2),\,d\}}+a(x))&\geq&-
d\epsilon+\int_{-\infty}^{+\infty}a(x)\delta(x-x_0)dx\\[2mm]
&=&-d\epsilon+a(x_0)\\[2mm]
&\geq&-d\epsilon+\sup\limits_{x\in\mathbb{R}}a(x)-\epsilon.
\end{array}
$$
It follows from the arbitrarily of $\epsilon$ that
$$
\liminf\limits_{-L_1, L_2\to+\infty}\lambda_{p}(\mathcal{L}_{\{(L_1,\,L_2),\,d\}}+a(x))\geq\sup\limits_{x\in\mathbb{R}}a(x).
$$

\epf

In the following, we discuss the monotonicity of the principal eigenvalue with respect to $d$ and its  limiting behaviors as $d\rightarrow0$ or $d\rightarrow+\infty$.

\begin{thm}
\label{thm3.3}
Suppose that $\mathbf{(J)}$ and $\mathbf{(H)}$ hold.
Then, the following statements hold:

$(i)$ $\lambda_{p}(\mathcal{L}_{\{(L_1,\,L_2),\,d\}}+a(x))$ is a strictly monotone decreasing function of $d$.

$(ii)$ $\lim\limits_{d\to 0}\lambda_{p}(\mathcal{L}_{\{(L_1,\,L_2),\,d\}}+a(x))=\max\limits_{x\in[L_1,\,L_2]}a(x)$.

$(iii)$ If $\int_{-\infty}^{L_1-L_2}J(z)dz>0$(symmetrically, $\int_{L_2-L_1}^{+\infty}J(z)dz>0$) holds,
then $\lim\limits_{d\to+\infty}\lambda_{p}(\mathcal{L}_{\{(L_1,\,L_2),\,d\}}+a(x))=-\infty$.
\end{thm}

\bpf
$(i)$ Assume that $\lambda_{p}(d_1):=\lambda_{p}(\mathcal{L}_{\{(L_1,L_2),\,d\}}+a(x))$ is the principal eigenvalue and $\phi(x)$
is the corresponding positive eigenfunction with $\|\phi\|_{L^{2}}=1$, we have
$$
\lambda_{p}(d_{1})\phi(x)=d_{1}\int_{L_1}^{L_2}J(x-y)\phi(y)dy
-d_{1}\phi(x)+a(x)\phi(x),\, \, \, x\in(L_1,L_2).
$$
Suppose that $d<d_1$, then
$$
\begin{array}{llllll}
\lambda_{p}(d_1)&=&d_1\int_{L_1}^{L_2}\int_{L_1}^{L_2}J(x-y)\phi(y)\phi(x)dydx
-d_1+\int_{L_1}^{L_2}a(x)\phi^2(x)dx\\[2mm]
&=&d\int_{L_1}^{L_2}\int_{L_1}^{L_2}J(x-y)\phi(y)\phi(x)dydx-d_1+
\int_{L_1}^{L_2}a(x)\phi^2(x)dx\\[2mm]
&&+(d_1-d)\int_{L_1}^{L_2}\int_{L_1}^{L_2}J(x-y)\phi(y)\phi(x)dydx\\[2mm]
&<&d\int_{L_1}^{L_2}\int_{L_1}^{L_2}J(x-y)\phi(y)\phi(x)dydx-d_1
+\int_{L_1}^{L_2}a(x)\phi^2(x)dx+d_1-d\\[2mm]
&=&d\int_{L_1}^{L_2}\int_{L_1}^{L_2}J(x-y)\phi(y)\phi(x)dydx
-d+\int_{L_1}^{L_2}a(x)\phi^2(x)dx\\[2mm]
&\leq&\lambda_{p}(d).
\end{array}
$$
Therefore $\lambda_{p}(d_1)<\lambda_{p}(d)$.

$(ii)$ The idea of this proof is from Theorem 2.8 in \cite{YLR}. For the eigenvalue problem
\begin{equation}
d\int_{L_1}^{L_2}J(x-y)\varphi(y)dy-d\varphi(x)+(\max\limits_{x\in[L_1,\,L_2]} a(x))\,\varphi(x)=\lambda^*_{p}\varphi(x),
\,\,x\in(L_1,L_2).
\label{th3.3-1}
\end{equation}
It follows from \cite{GR} that $\lambda^*_{p}\leq\max\limits_{x\in[L_1,\,L_2]} a(x)$. Therefore,
we have $\lambda_{p}(\mathcal{L}_{\{(L_1,\,L_2),\,d\}}+a(x))\leq\lambda^*_{p}\leq\max\limits_{x\in[L_1,\,L_2]} a(x)$ by $(ii)$ of Proposition \ref{prop3.1}.

Next, we prove that $\liminf\limits_{d\to 0}\lambda_{p}(\mathcal{L}_{\{(L_1,\,L_2),\,d\}}+a(x))\geq\max\limits_{x\in[L_1,\,L_2]} a(x)$.
Assume for the contrary that $\liminf\limits_{d\to 0}\lambda_{p}(\mathcal{L}_{\{(L_1,\,L_2),\,d\}}+a(x))\leq\max\limits_{x\in[L_1,\,L_2]} a(x)-\epsilon$
for some $\epsilon>0$. By the definition of $\liminf$, there exists some $\widehat{d}>0$ such that if $d\leq \widehat{d}$, then
$$
\lambda_{p}(\mathcal{L}_{\{(L_1,\,L_2),\,d\}}+a(x))\leq\max\limits_{x\in[L_1,\,L_2]} a(x)-\frac{\epsilon}{2}.
$$
On the other hand, by the continuity of $a(x)$, there exist $x_0\in(L_1,L_2)$ and $r>0$ such that
$$
\max\limits_{x\in[L_1,\,L_2]} a(x)\leq a(x)+\frac{\epsilon}{4}, \, \, x\in U_{r}(x_0)\subset(L_1,L_2).
$$
Therefore,
$$
\lambda_{p}(\mathcal{L}_{\{(L_1,\,L_2),\,d\}}+a(x))\leq a(x)-\frac{\epsilon}{4}
$$
for $0<d<\widehat{d}$ and $x\in U_{r}(x_0)$. Let $(\lambda_{p}(\mathcal{L}_{\{(L_1,\,L_2),\,d\}}+a(x)),\psi(x))$ be the eigenpair of the following eigenvalue problem:
$$
d\int_{L_1}^{L_2}J(x-y)\psi(y)dy-d\psi(x)+ a(x)\psi(x)=\lambda_{p}(\mathcal{L}_{\{(L_1,\,L_2),\,d\}}+a(x))\psi(x),\ x\in(L_1,L_2).
$$
Then,
$$
\int_{L_1}^{L_2}J(x-y)\psi(y)dy-\psi(x)=\frac{\lambda_{p}(\mathcal{L}_{\{(L_1,\,L_2),\,d\}}
+a(x))-a(x)}{d}\psi(x)\leq-\frac{\epsilon}{4d}\psi(x) \,\, \, \mbox{in}\, \, \, U_{r}(x_0).
$$
Let $\widetilde{\lambda}$ be the principal eigenvalue of the linear problem
\begin{equation}
\left\{\begin{array}{lll}
\int_\mathbb{R}J(x-y)u(y)dy-u(x)=\lambda u(x)\,\;&\,\mbox{in}\, \, \, U_{r}(x_0), \\[2mm]
u(x)=0\;&\, \mbox{in}\, \, \, \mathbb{R}\backslash U_{r}(x_0).
\end{array} \right.
\label{th3.3-2}
\end{equation}
It is well known that $-1<\widetilde{\lambda}<0$ by Theorem 2.1 in \cite{GR1}. Let $\Psi(x)$ be
the eigenfunction corresponding to $\widetilde{\lambda}$ and $\|\Psi(x)\|_{L^{\infty}}=1$. Take
$$
\overline{\psi}(x)=\frac{\psi(x)}{\inf_{U_{r}(x_0)}\psi(x)},\, \, \, \underline{\psi}(x)=\Psi(x).
$$
We consider the following problem
\begin{equation}
\left\{\begin{array}{lll}
\int_\mathbb{R}J(x-y)u(y)dy-u(x)=-\frac{\epsilon}{4d} u(x)\,\;&\,\mbox{in}\, \, \, U_{r}(x_0), \\[2mm]
u(x)=0\;&\, \mbox{in}\, \, \, \mathbb{R}\backslash U_{r}(x_0).
\end{array} \right.
\label{th3.3-3}
\end{equation}
Direct calculation yields
$$
\begin{array}{llllll}
&&\int_{U_{r}(x_0)}J(x-y)\overline{\psi}(y)dy-\overline{\psi}(x)+\frac{\epsilon}{4d}\overline{\psi}(x)\\[2mm]
&=&\frac{1}{\inf\psi}[\int_{U_{r}(x_0)}J(x-y)\psi(y)dy-\psi(x)]+
\frac{\epsilon}{4d}\overline{\psi}(x)\\[2mm]
&\leq&\frac{1}{\inf\psi}(-\frac{\epsilon}{4d}\psi(x)+\frac{\epsilon}{4d}\psi(x))\\[2mm]
&=&0,
\end{array}
$$
and
$$
\begin{array}{llllll}
&&\int_{U_{r}(x_0)}J(x-y)\underline{\psi}(y)dy-\underline{\psi}(x)+\frac{\epsilon}{4d}\underline{\psi}(x)\\[2mm]
&=&\int_{U_{r}(x_0)}J(x-y)\Psi(x)dy-\Psi(x)]+
\frac{\epsilon}{4d}\Psi(x)\\[2mm]
&=&\widetilde{\lambda}\Psi(x)+\frac{\epsilon}{4d}\Psi(x)\\[2mm]
&\geq&0
\end{array}
$$
provided $d<\min\{\widehat{d},-\frac{\epsilon}{4\widetilde{\lambda}}\}$. Hence, by the super-sub solution method in \cite{GR2}, one can yield \eqref{th3.3-3} has a positive solution between $\overline{\psi}(x)$ and $\underline{\psi}$, which implies that $\widetilde{\lambda}=-\frac{\epsilon}{4d}$. This contradicts to the independence of $\widetilde{\lambda}$ from $d$. Therefore, $\lim\limits_{d\to 0}\lambda_{p}(\mathcal{L}_{\{(L_1,\,L_2),\,d\}}+a(x))=\max\limits_{x\in[L_1,\,L_2]}a(x)$.

$(iii)$ We claim that if $\int_{-\infty}^{L_1-L_2}J(z)dz>0$ (or $\int_{L_2-L_1}^{+\infty}J(z)dz>0$), then  $\lim\limits_{d\to+\infty}\lambda_{p}(\mathcal{L}_{\{(L_1,\,L_2),\,d\}}+a(x)):=\lambda_{\infty}=-\infty$. We argue by contradiction and suppose that $\lambda_{\infty}>-\infty$. Now let $\varphi(x)=C_1$ (positive constant), without loss of generality, taking $C_1=1$,
then for $x\in(L_1,L_2)$,
$$
\begin{array}{llllll}
&&d\int_{L_1}^{L_2}J(x-y)\varphi(y)dy
-d\varphi(x)+a(x)\varphi(x)\\[2mm]
&<&d\int_{L_1}^{L_2}J(x-y)dy-d+\max\limits_{x\in[L_1,L_2]}a(x)\\[2mm]
&=&-d\int_{\mathbb{R}\backslash[L_1,L_2]}J(x-y)dy+\max\limits_{x\in[L_1,L_2]}a(x)\\[2mm]
&=&-d(\int_{-\infty}^{x-L_2}+\int_{x-L_1}^{+\infty})J(z)dz+\max\limits_{x\in[L_1,L_2]}a(x)\\
&\leq&-d(\int_{-\infty}^{L_1-L_2}+\int_{L_2-L_1}^{+\infty})J(z)dz+\max\limits_{x\in[L_1,L_2]}a(x)
\end{array}
$$
as $\int_{\mathbb{R}}J(x)dx=1$. Owing to the assumption that $\int_{-\infty}^{L_1-L_2}J(z)dz>0$ (or $\int_{L_2-L_1}^{+\infty}J(z)dz>0$), there exists a $d$ adequately large such that
$$
d\int_{L_1}^{L_2}J(x-y)\varphi(y)dy
-d\varphi(x)+a(x)\varphi(x)\leq(\lambda_{\infty}-1)\varphi(x).
$$
Thus, by definition of $\lambda_{p}(\mathcal{L}_{\{(L_1,\,L_2),\,d\}}+a(x))$, we have
$$
\lambda_{p}(\mathcal{L}_{\{(L_1,\,L_2),\,d\}}+a(x))\leq\lambda_{\infty}-1,
$$
further,
$$
\lim\limits_{d\to+\infty}\lambda_{p}(\mathcal{L}_{\{(L_1,\,L_2),\,d\}}+a(x))
=\lambda_{\infty}\leq\lambda_{\infty}-1.
$$
We obtain the desired contradiction. Hence, $\lim\limits_{d\to+\infty}\lambda_{p}(\mathcal{L}_{\{(L_1,\,L_2),\,d\}}+a(x))=-\infty$.

\epf

\begin{rmk}
\label{rmk3.1}
Compared with \cite{YLR}, where the nonlocal operator is $d\int_{\Omega}J(x-y)(\psi(y)-\psi(x))dy$, we
consider the nonlocal operator $d\int_{\Omega}J(x-y)\varphi(y)dy-d\varphi(x)$ and
prove that when $d\rightarrow+\infty$, the limit of the principal eigenvalue is $-\infty$ for
some cases, which is different from the result in \cite{YLR}; its limit is the
average of $a(x)$ over $\Omega$.
\end{rmk}

\section{Spreading-vanishing }

Since $h(t)$ and $-g(t)$ are monotonically increasing with $t>0$, there exist $h_{\infty}$
and $g_{\infty}$ such that $\lim\limits_{t\rightarrow+\infty}g(t)=g_{\infty}\in[-\infty,-h_{0})$
and $\lim\limits_{t\rightarrow+\infty}h(t)=h_{\infty}\in(h_{0},+\infty]$. Here we define that  $\mathbf{vanishing}$ occurs if $h_{\infty}-g_{\infty}<+\infty$ and $\lim\limits_{t\rightarrow+\infty}\max\limits_{x\in[g(t),\,h(t)]}I(t,x)=0$;
and $\mathbf{spreading}$ happens provided that $h_{\infty}-g_{\infty}=+\infty$ and $\limsup\limits_{t\rightarrow+\infty}\|I(\cdot,t)\|_{C([g(t),\,h(t)])}>0$.
In this section, we always assume that $\mathbf{(J)}$ and $\mathbf{(H)}$ hold. The following proposition directly holds from Theorem 3.5 in \cite{CLW}.

\begin{prop}
\label{prop5.1}
Let $(S,I;g,h)$ be the unique solution of \eqref{a02}. If $h_{\infty}-g_{\infty}<+\infty$, then
$\lim\limits_{t\rightarrow+\infty}g'(t)=\lim\limits_{t\rightarrow+\infty}h'(t)=0$.
\end{prop}

Next, we discuss the asymptotic behavior of the solution of problem \eqref{a02} when
$h_{\infty}-g_{\infty}<+\infty$.

\begin{thm}
\label{thm5.2}
Let $(S,I;g,h)$ be the unique solution of problem \eqref{a02} with $h_{\infty}-g_{\infty}<+\infty$,
then $\lim\limits_{t\rightarrow+\infty}\max\limits_{x\in[g(t),\,h(t)]}I(t,x)=0$, $\lim\limits_{t\rightarrow+\infty}S(t,x)=\frac{\sigma}{\mu_1}$ and $\lambda_{p}(\mathcal{L}_{\{(g_{\infty},\,h_{\infty}),\,d\}}+a(x))\leq 0$.
\end{thm}

\bpf
Assume by contradiction that $\lim\limits_{t\rightarrow+\infty}\max\limits_{x\in[g(t),\,h(t)]}I(t,x)>0$, there exists $\epsilon_{1}>0$ and sequence $\{(t_i,x_i)\}_{i=1}^{\infty}$ with $x_{i}\in[g(t),h(t)]$ and $t_{i}\rightarrow+\infty$ as $i\rightarrow+\infty$ such that
$I(t_{i},x_{i})\geq\frac{\epsilon_{1}}{2}$ for $i\in\mathbb{N}$. Since $g_{\infty}<g(t)<x_{i}<h(t)<h_{\infty}$, there exists a subsequence $\{x_{i_{j}}\}_{j=1}^{\infty}$
such that $x_{i_{j}}\rightarrow x_0\in(g_{\infty},h_{\infty})$ as $j\rightarrow +\infty$. For $t\in(-t_{i},+\infty)$ and $x\in(g(t+t_i),h(t+t_i))$, define
$$
\overline{I_{i}}(t,x)=I(t+t_{i},x).
$$
Applying Theorem \ref{thm2.2} gives that $I$ and $S$ are positive and bounded, and then $\overline{I_{i}}(t,x)$ satisfies
$$
\overline{I_{it}}\geq d\int_{g_{i}(t)}^{h_{i}(t)}J(x-y)\overline{I_{i}}(t,y)dy-d \overline{I_{i}}(t,x)-\mu_{2}\overline{I_{i}}-\gamma(b,0,x)\overline{I_{i}},\,
t>-t_i,\, x\in(g_{i}(t), h_{i}(t)).
$$
We next consider the following auxiliary problem
{\small
\begin{equation*}
\left\{\begin{array}{lll}
u_{t}=d\int_{g_{i}(t)}^{h_{i}(t)}J(x-y)u(t,y)dy-d u(t,x)-\mu_{2}u-\gamma(b,0,x)u,\,\;
&\,t>-t_i,\, x\in(g_{i}(t), h_{i}(t)),\\[2mm]
u(0,x)=\overline{I}_{i}(0,x),\,\;&\,x\in(g_{i}(t), h_{i}(t)),
\end{array} \right.
\end{equation*}}
it follows that $u(t,x)\rightarrow U(t,x)$ as $i\rightarrow+\infty$, and $U(t,x)$ satisfies
{\small
\begin{equation*}
\left\{\begin{array}{lll}
U_{t}(t,x)=d\int_{g_\infty}^{h_\infty}J(x-y)U(t,y)dy-dU(t,x)-\mu_{2}U-\gamma(b,0,x)U,\,\;
&\,t\in\mathbb{R},\, \, \, x\in(g_\infty,h_\infty), \\[2mm]
U(0,x_0)=\lim\limits_{i\rightarrow+\infty}\overline{I}_{i}(0,x_{i})=
\lim\limits_{i\rightarrow+\infty}I(t_{i},x_{i})\geq\frac{\epsilon_{1}}{2}>0,
\end{array} \right.
\end{equation*}}
and then $U(t,x)>0$ in $\mathbb{R}\times(g_\infty,h_\infty)$ by the maximum principle \cite{DN} for the nonlocal problem.

On the other hand, considering $h_{\infty}-g_{\infty}<+\infty$ and Proposition \ref{prop5.1}, we have  $\lim\limits_{t\rightarrow+\infty}g'(t)=\lim\limits_{t\rightarrow+\infty}h'(t)=0$ as $t\rightarrow+\infty$, which means
$$
\begin{array}{llllll}
0=\lim\limits_{i\rightarrow+\infty}h'(t+t_{i})&=&
k\lim\limits_{i\rightarrow+\infty}\int_{g(t+t_{i})}^{h(t+t_{i})}
\int_{h(t+t_{i})}^{+\infty}J(x-y)\overline{I}_{i}(t,x)dydx
\\[2mm]
&\geq&k\int_{g(\infty)}^{h(\infty)}
\int_{h(\infty)}^{+\infty}J(x-y)U(t,x)dydx\\[2mm]
&>&0
\end{array}
$$
and
$$
\begin{array}{llllll}
0=\lim\limits_{i\rightarrow+\infty}g'(t+t_{i})&=&
-k\lim\limits_{i\rightarrow+\infty}\int_{g(t+t_{i})}^{h(t+t_{i})}
\int_{-\infty}^{g(t+t_{i})}J(x-y)\overline{I}_{i}(t,x)dydx\\[2mm]
&\leq&-k\int_{g(\infty)}^{h(\infty)}
\int_{-\infty}^{g(\infty)}J(x-y)U(t,x)dydx\\[2mm]
&<&0.
\end{array}
$$
It is a contradiction. Hence, $\lim\limits_{t\rightarrow+\infty}\max\limits_{x\in[g(t),\,h(t)]}I(t,x)=0$.

Next, we will prove that $\lim\limits_{t\rightarrow+\infty}S(t,x)=\frac{\sigma}{\mu_{1}}$. Since $\lim\limits_{t\rightarrow+\infty}\max\limits_{x\in[g(t),\,h(t)]}I(t,x)=0$, for any $\epsilon>0$, one can choose a $T>0$ large such that
$$
0<I(t,x)<\epsilon
$$
for $t>T$, $x\in(g(t),h(t))$.

Obviously, $S(t,x)$ satisfies
$$
S_{t}\geq d\mathcal{L}_1[S]+\sigma-\mu_{1}S-\epsilon\beta(m,0,x)S,\, \, t>T,\ x\in \mathbb{R},
$$
and then $S(t,x)\geq\underline{S}(t)$, where $\underline{S}(t)$ is the solution to problem
\begin{equation*}
\left\{ \begin{array}{ll}
\underline{S}_{t}=\sigma-(\mu_{1}+\epsilon\sup\limits_{x\in\mathbb{R}}\beta(m,0,x))\underline{S}, &t>T,\\[2mm]
\underline{S}(T)=\inf\limits_{x\in\mathbb{R}}S(T,x).
\end{array}\right.
\end{equation*}
It follows from Lemma 2.4 in \cite{LSW} that $\lim\limits_{t\rightarrow+\infty}\underline{S}(t)=\frac{\sigma}{\mu_{1}+
\epsilon\sup\limits_{x\in\mathbb{R}}\beta(m,0,x)}$.
Therefore, $\liminf\limits_{t\rightarrow+\infty}S(t,x)\geq\frac{\sigma}{\mu_{1}+
\epsilon\sup\limits_{x\in\mathbb{R}}\beta(m,0,x)}$. Letting $\epsilon\rightarrow0$ yields
\begin{equation}
\liminf\limits_{t\rightarrow+\infty}S(t,x)\geq\frac{\sigma}{\mu_{1}}.
\label{th5.2-1}
\end{equation}
On the other hand, $S(t,x)$ satisfies
$$
S_{t}\leq d\mathcal{L}_1[S]+\sigma-\mu_{1}S+\gamma(b,\epsilon,x)\epsilon,\, \, t>T,\ x\in \mathbb{R}.
$$
Let $\overline{S}(t)$ be the solution of
\begin{equation*}
\left\{ \begin{array}{ll}
\overline{S}_{t}=\sigma-\mu_{1}\overline{S}(t)+\sup\limits_{x\in\mathbb{R}}\gamma(b,\epsilon,x)\epsilon, &t>T,\\[2mm]
\overline{S}(T)=\sup\limits_{x\in\mathbb{R}}S(T,x).
\end{array}\right.
\end{equation*}
Apparently, $S(t,x)\leq\overline{S}(t)$ by Lemma 2.4 in \cite{LSW} and $\lim\limits_{t\rightarrow+\infty}\overline{S}(t)=\frac 1{\mu_{1}}[\sigma+
\sup\limits_{x\in\mathbb{R}}\gamma(b,\epsilon,x)\epsilon]$. Therefore,
$$
\limsup\limits_{t\rightarrow+\infty}S(t,x)\leq\lim\limits_{t\rightarrow+\infty}\overline{S}(t)=\frac 1{\mu_{1}}[\sigma+
\sup\limits_{x\in\mathbb{R}}\gamma(b,\epsilon,x)\epsilon].
$$
Letting $\epsilon\rightarrow0$ gives
\begin{equation}
\limsup\limits_{t\rightarrow+\infty}S(t,x)\leq\frac{\sigma}{\mu_{1}},
\label{th5.2-2}
\end{equation}
which together with \eqref{th5.2-1} yields $\lim\limits_{t\rightarrow+\infty}S(t,x)=\frac{\sigma}{\mu_{1}}$ uniformly for $x\in \mathbb{R}$.

In what follows, we prove that $\lambda_{p}(\mathcal{L}_{\{(g_{\infty},\,h_{\infty}),\,d\}}+a(x))\leq 0$. We argue by contradiction and suppose that $\lambda_{p}(\mathcal{L}_{\{(g_{\infty},\,h_{\infty}),\,d\}}+a(x))>0$. Owing to the continuous dependence of $\lambda_{p}(\mathcal{L}_{\{(g_{\infty},\,h_{\infty}),\,d\}}+a(x))$ on $a(x)$ and $(g_{\infty},h_{\infty})$, there exists a small $\epsilon$ such that $\lambda_{p}(\mathcal{L}_{\{(g_{\infty}+\epsilon,\,h_{\infty}-\epsilon),\,d\}}+a(x)-\beta(m,0,x)\epsilon)>0$. Furthermore, in view of $S(t,x)\rightarrow\frac{\sigma}{\mu_{1}}$ as $t\rightarrow+\infty$ and $h_{\infty}-g_{\infty}<+\infty$, there exists a $T^*>0$ such that
$$
g(t)<g_{\infty}+\epsilon,\,\,h(t)>h_{\infty}-\epsilon,\,\, t>T^*,
$$
$$
S(t,x)>\frac{\sigma}{\mu_{1}}-\epsilon,\, \, t>T^*,\,\,x\in\mathbb{R}.
$$
Then, for $t>T^*$, $x\in(g(t),h(t))$,
$$
\begin{array}{llllll}
I_{t}(t,x)&=&d\int_{g(t)}^{h(t)}J(x-y)I(t,y)dy-d I(t,x)-\mu_{2}I+\beta(m,I,x)SI-\gamma(b,I,x)I\\[2mm]
&\geq&d\int_{g_{\infty}+\epsilon}^{h_{\infty}-\epsilon}J(x-y)I(t,y)dy-d I(t,x)-\mu_{2}I+\beta(m,I,x)(\frac{\sigma}{\mu_{1}}-\epsilon)I-\gamma(b,I,x)I\\[2mm]
&\geq&d\int_{g_{\infty}+\epsilon}^{h_{\infty}-\epsilon}J(x-y)I(t,y)dy-d I(t,x)-\mu_{2}I+\beta(m,0,x)(\frac{\sigma}{\mu_{1}}-\epsilon)I-\gamma(b,0,x)I.
\end{array}
$$
Let $\phi(x)$ be the eigenfunction corresponding to $\lambda_{p}(\mathcal{L}_{\{(g_{\infty}+\epsilon,\,h_{\infty}-\epsilon),\,d\}}+a(x)-\beta(m,0,x)\epsilon)$ and $\|\phi(x)\|_{L^{\infty}}=1$, and then $\phi(x)$ satisfies
$$
\begin{array}{llllll}
&&d\int_{g_{\infty}+\epsilon}^{h_{\infty}-\epsilon}J(x-y)\phi(y)dy-d \phi-\mu_{2}\phi+\beta(m,0,x)(\frac{\sigma}{\mu_{1}}-\epsilon)\phi-\gamma(b,0,x)\phi\\[2mm]
&=&d\int_{g_{\infty}+\epsilon}^{h_{\infty}-\epsilon}J(x-y)\phi(y)dy-d\phi+a(x)\phi
-\beta(m,0,x)\epsilon\phi\\[2mm]
&=&\lambda_{p}(\mathcal{L}_{\{(g_{\infty}+\epsilon,\,h_{\infty}-\epsilon),\,d\}}
+a(x)-\beta(m,0,x)\epsilon)\phi.
\end{array}
$$
If we choose $\delta$ sufficiently small such that $\delta\phi(x)\leq I(T^*,x)$ for $x\in[g_{\infty}+\epsilon,h_{\infty}-\epsilon]$, then
$$
I(t,x)\geq\delta\phi(x)>0\, \, \, \mbox{for}\, \, \, t>T^*,\, \, \,x\in[g_{\infty}+\epsilon,h_{\infty}-\epsilon]
$$
by the comparison principle in \cite{DN}, which leads to a contradiction to the fact $\lim\limits_{t\rightarrow+\infty}\max\limits_{x\in[g(t),\,h(t)]}I(t,x)=0$.

\epf

\begin{thm}
\label{thm5.3}
Suppose $\lambda_{p}(\mathcal{L}_{\{(-h_{0},\,h_{0}),\,d\}}+a(x))<0$. Then $h_{\infty}-g_{\infty}<+\infty$,
$\lim\limits_{t\rightarrow+\infty}\max\limits_{x\in[g(t),\,h(t)]}I(t,x)=0$ and $\lim\limits_{t\rightarrow+\infty}S(t,x)=\frac{\sigma}{\mu_1}$ if $\|S_0(x)\|_{L^{\infty}(\mathbb{R})}+\|I_{0}(x)\|_{C([-h_0,\,h_0])}$
is sufficiently small.
\end{thm}

\bpf
Applying Lemma \ref{lem2.1} yields $S(t,x)\leq\frac{\sigma}{\mu_1}$ for $t>0,\,x\in\mathbb{R}$ if $\|S_0(x)\|_{L^{\infty}(\mathbb{R})}+\|I_{0}(x)\|_{C([-h_0,\,h_0])}\leq\frac{\sigma}{\mu_1}$.
In view of $\lambda_{p}(\mathcal{L}_{\{(-h_{0},\,h_{0}),\,d\}}+a(x))<0$, there exists a $\epsilon>0$ small such that
$\lambda_{p}(\mathcal{L}_{\{(-h_{\epsilon},\,h_{\epsilon}),\,d\}}+a(x))<0$ with $h_{\epsilon}=h_{0}+\epsilon$ by Theorem \ref{thm3.2} $(i)$. Let $\phi(x)$ be the eigenfunction of the principal eigenvalue $\lambda_{p}(\mathcal{L}_{\{(-h_{\epsilon},\,h_{\epsilon}),\,d\}}+a(x))$, which satisfies
$$
d\int_{-h_{\epsilon}}^{h_{\epsilon}}J(x-y)\phi(y)dy-d\phi(x)+a(x)\phi(x)=\lambda_{p}
(\mathcal{L}_{\{(-h_{\epsilon},\,h_{\epsilon}),\,d\}}+a(x))\phi(x),
\,\,x\in(-h_{\epsilon},h_{\epsilon}).
$$
Denote
$$
M=\delta C(\int_{-h_{\epsilon}}^{h_{\epsilon}}\phi(x)dx)^{-1}, \, \, \, C=\frac{h_{\epsilon}-h_{0}}{k},
$$
$$
\overline{h}(t)=h_{0}+kC[1-e^{-\delta t}], \, \, \, \overline{g}(t)=-\overline{h}(t),
$$
$$
\overline{h}'(t)=kC\delta e^{-\delta t}, \, \, \, \overline{I}(t,x)=Me^{-\delta t}\phi(x).
$$
By direct calculations, we have
$$
\begin{array}{llllll}
&&k\int_{\overline{g}}^{\overline{h}}\int_{\overline{h}}^{\infty}J(x-y)\overline{I}(t,x)dydx\\[2mm]
&\leq&k\int_{\overline{g}}^{\overline{h}}\overline{I}(t,x)dx\\[2mm]
&=&k\delta Ce^{-\delta t}=\overline{h}'(t),\, \, \, \, t>0.
\end{array}
$$
In a similar way, one can deduce $-k\int_{\overline{g}}^{\overline{h}}\int_{-\infty}^{\overline{g}}J(x-y)\overline{I}(t,x)dydx\geq
\overline{g}'(t)$ for $t>0$. Clearly, $\overline{I}(t,\overline{g}(t))>0$ and $\overline{I}(t,\overline{h}(t))>0$ for $t>0$. For $t>0$ and $x\in(\overline{g},\overline{h})$, we obtain
$$
\begin{array}{llllll}
&&\overline{I}_{t}(t,x)-d\int_{\overline{g}(t)}^{\overline{h}(t)}J(x-y)\overline{I}(t,y)dy+d \overline{I}(t,x)+\mu_{2}\overline{I}-\beta(m,\overline{I},x)\overline{S}\overline{I}+
\gamma(b,\overline{I},x)\overline{I}\\[2mm]
&\geq&-[d\int_{\overline{g}(t)}^{\overline{h}(t)}J(x-y)Me^{-\delta t}\phi(y)dy-d Me^{-\delta t}\phi(x)-\mu_{2}Me^{-\delta t}\phi(x)]\\[2mm]
&&-\delta\overline{I}+\gamma(b,\overline{I},x)\overline{I}-\beta(m,\overline{I},x)\overline{S}\overline{I}\\[2mm]
&\geq&-\delta\overline{I}-(\lambda_{p}(\mathcal{L}_{(-h_{\epsilon},\,h_{\epsilon})}+a(x))+\gamma(b,0,x)-
\beta(m,0,x)\frac{\sigma}{\mu_1})\overline{I}+\gamma(b,\overline{I},x)\overline{I}-
\beta(m,\overline{I},x)\frac{\sigma}{\mu_1}\overline{I}\\[2mm]
&=&\overline{I}(-\delta-\lambda_{p}(\mathcal{L}_{(-h_{\epsilon},\,h_{\epsilon})}+a(x))-\gamma(b,0,x)+
\gamma(b,\overline{I},x)+\frac{\sigma}{\mu_1}(\beta(m,0,x)-\beta(m,\overline{I},x))).
\end{array}
$$
Recalling that $\gamma(b,\overline{I},x)\rightarrow\gamma(b,0,x)$ and $\beta(m,\overline{I},x)\rightarrow\beta(m,0,x)$ as $\delta\rightarrow 0$,
we can choose $\delta$ small enough so that
$$
\overline{I}_{t}(t,x)-d\int_{\overline{g}(t)}^{\overline{h}(t)}J(x-y)\overline{I}(t,y)dy+d \overline{I}(t,x)+\mu_{2}\overline{I}-\beta(x)\overline{S}\overline{I}+
\gamma(x,b,\overline{I})\overline{I}\geq0.
$$
Moreover, if $\|S_0(x)\|_{L^{\infty}(\mathbb{R})}+\|I_{0}(x)\|_{C([-h_0,\,h_0])}$ sufficiently small,
then $I_{0}\leq M\phi(x)$, $x\in[-h_{0},h_{0}]$.

Applying the comparison principle gives
$$
\overline{g}(t)\leq g(t),\, \, \, h(t)\leq\overline{h}(t)\, \rm{and}\, \, \ I(t,x)\leq\overline{I}(t,x)
$$
for $t>0$ and $g(t)<x<h(t)$. Therefore,
\begin{equation}
\lim\limits_{t\rightarrow+\infty}I(t,x)\leq\lim\limits_{t\rightarrow+\infty}\overline{I}(t,x)=0
\label{thm5}
\end{equation}
and
$$
h_{\infty}-g_{\infty}\leq2h_{\epsilon}<+\infty,
$$
which gives that $\lim\limits_{t\rightarrow+\infty}S(t,x)=\frac{\sigma}{\mu_1}$ uniformly for $x\in\mathbb{R}$ by Theorem \ref{thm5.2}.\\[0.5mm]
\epf

\begin{rmk}
\label{rmk1.1}
It follows from the proof of Theorem \ref{thm5.3} that $M\rightarrow+\infty$ as $k\rightarrow0$. Therefore, there exists a $k_{*}>0$ such that \eqref{thm5} holds for all $k\in(0,k_{*})$ for any given initial function pair $(S_{0}(x),I_{0}(x))$.
\end{rmk}

\begin{thm}
\label{thm5.4}
If $-g_\infty=h_\infty=+\infty$ and $\sup\limits_{x\in\mathbb{R}}a(x)>0$, then $\limsup\limits_{t\rightarrow+\infty}\|I(\cdot,t)\|_{C([g(t),\,h(t)])}>0$.
\end{thm}
\bpf Conversely, suppose that $\lim\limits_{t\rightarrow+\infty}\|I(t,\cdot)\|_{C([g(t),\,h(t)])}=0$.
It follows from Theorem \ref{thm5.2} that
\begin{equation}
\lim\limits_{t\rightarrow+\infty}S(t,x)=\frac{\sigma}{\mu_1} \, \, \ \textrm{uniformly for} \, \, \ x\in\mathbb{R}.
\label{d-21}
\end{equation}
Since $-g_\infty=h_\infty=+\infty$, from $(iii)$ of Theorem \ref{thm3.2}, there exists a $T^*>0$ large enough such that
$$
\lambda_{p}(\mathcal{L}_{\{(g(T^*),\,h(T^*)),\,d\}}+a(x))>0\, \, \, \ \textrm{for}\, \, \, \ t\geq T^*.
$$
Now, we consider the eigenvalue problem
\small{$$
d\int_{g(T^*)}^{h(T^*)}J(x-y)\psi(y)dy-d\psi(x)+a(x)\psi(x)=
\lambda_{p}(\mathcal{L}_{\{(g(T^*),\,h(T^*)),\,d\}}+a(x))\psi(x),\,\,x\in(g(T^*),h(T^*)),
$$}
and positive function $\psi(x)$ with $\|\psi\|_{L^\infty((g(T^*),\,h(T^*)))}=1$ is its eigenfunction to
the principal eigenvalue $\lambda_{p}(\mathcal{L}_{\{(g(T^*),\,h(T^*)),\,d\}}+a(x))$.

In view of \eqref{d-21}, for any given $0<\epsilon<\min\{\frac{\sigma}{\mu_1},\,\lambda_{p}(\mathcal{L}_{\{(g(T^*),\,h(T^*)),\,d\}}+a(x))
(\sup\limits_{x\in\mathbb{R}}\beta)^{-1}\}$, there exists a $T^{**}>T^*$ such
that
$$
S(t,x)>\frac{\sigma}{\mu_1}-\epsilon,\, \, \, \ t\geq T^{**}, \, \, \, \ x\in[g(T^*),h(T^*)],
$$
and then $I(t,x)$ satisfies
\begin{equation*}
\left\{\begin{array}{lll}
I_{t}\geq d\int_{g(T^*)}^{h(T^*)}J(x-y)I(t,y)dy-d I(t,x)-\mu_{2}I\\[2mm]
\qquad +\beta(m,I,x)(\frac{\sigma}{\mu_1}-\epsilon)I-\gamma(b,I,x)I,\; &\, t>T^{**},\,\, x\in(g(T^*),h(T^*)), \\[2mm]
%I(g(T^*),t)>0, \, \ I(h(T^*),t)>0,\;&\, t>T^{**}, \\[2mm]
I(T^{**},x)>0,\;&\, g(T^*)\leq x\leq h(T^*).
\end{array} \right.
\end{equation*}
We now construct a suitable lower solution for the following auxiliary problem
\begin{equation}
\left\{\begin{array}{lll}
W_{t}=d\int_{g(T^*)}^{h(T^*)}J(x-y)W(t,y)dy-d W(t,x)\\[2mm]
\qquad -\mu_{2}W+\beta(m,W,x)(\frac{\sigma}{\mu_1}-\epsilon)W-\gamma(b,W,x)W,\; &\, t>T^{**},\ x\in(g(T^*),h(T^*)), \\[2mm]
W(T^{**},x)=I(T^{**},x),\;&\, g(T^*)\leq x\leq h(T^*).
\end{array} \right.
\label{*}
\end{equation}
Choose
$$
\underline{W}(t,x)=\delta\psi(x),\, \, \, \, \ t>T^{**},\, \, \ g(T^*)\leq x\leq h(T^*),
$$
where $\delta>0$ is small enough such that $\delta\psi(x)\leq I(T^{**},x)$ for $x\in[g(T^*),h(T^*)]$.

For $t>T^{**}$ and $g(T^*)<x<h(T^*)$, direct computation yields
$$
\begin{array}{llllll}
&&\underline{W}_{t}-d\int_{g(T^*)}^{h(T^*)}J(x-y)\underline{W}(y)dy+d \underline{W}+\mu_{2}\underline{W}-\beta(m,\underline{W},x)(\frac{\sigma}{\mu_1}-\epsilon)\underline{W}
+\gamma(b,\underline{W},x)\underline{W}\\[2mm]
&=&-\delta(d\int_{g(T^*)}^{h(T^*)}J(x-y)\psi(t,y)dy-d \psi(x)-\mu_{2}\psi+\beta(m,\delta\psi,x)(\frac{\sigma}{\mu_1}-\epsilon)\psi-
\gamma(b,\delta\psi,x)\psi)\\[2mm]
&=&-\delta\psi(\lambda_{p}+\gamma(b,0,x)-
\beta(m,0,x)\frac{\sigma}{\mu_1}+\beta(m,\delta\psi,x)(\frac{\sigma}{\mu_1}-\epsilon)-
\gamma(b,\delta\psi,x))\\[2mm]
&\leq&\delta\psi(-\lambda_{p}(\mathcal{L}_{\{(g(T^*),h(T^*)),\,d\}}+a(x))+
\beta(m,\delta\psi,x)\epsilon)\\[4mm]
&<& 0.
\end{array}
$$
Note that the boundaries $g(T^*)$ and $h(T^*)$ are fixed, so there is no need to compare the boundary values of $W(t,x)$ and $I(t,x)$ by Lemma 3.1 in \cite{DN}.
Applying the comparison principle in \cite{DN} gives
$$
I(t,x)\geq W(t,x)\geq \underline{W}(t,x)=\delta\psi(x) \, \, \ \textrm{in}\, \, \ [T^{**},+\infty)\times(g(T^*),h(T^*)).
$$
Therefore,
$\liminf\limits_{t\rightarrow+\infty}I(t,x)\geq\liminf\limits_{t\rightarrow+\infty}W(t,x)\geq\delta\psi(0)>0$,
which is a contradiction.
\epf
\begin{rmk}
\label{rmk1.2}
Suppose that $a(x)$ is a positive constant. If $h_\infty-g_\infty=+\infty$, then $\limsup\limits_{t\rightarrow+\infty}\|I(\cdot,t)\|_{C([g(t),\,h(t)])}>0$. In the special case, similar
to Theorem 3.10 in \cite{CDLL}, a spreading-vanishing dichotomy holds.
\end{rmk}

Suppose that
$$
a(0)\geq d,
$$
we know that $\lambda_{p}(\mathcal{L}_{\{(L_1,\,L_2),\,d\}}+a(x))>0$ for any interval $(L_1,L_2)$ by  Theorem \ref{thm3.2} and $(i)$ of Proposition \ref{prop3.1}, which yields the following conclusion by Theorem \ref{thm5.2}.

\begin{thm}
\label{thm5.6}
If $a(0)\geq d$ holds, then spreading always occurs for \eqref{a02}.
\end{thm}

Next we consider the case $0<a(0)<d$. According to Theorem \ref{thm3.2} and $(i)$ of Proposition \ref{prop3.1}, there exists $L^*>0$ such that
\bess
\lambda_{p}(\mathcal{L}_{\{(L_1,\,L_2),\,d\}}+a(x))\left\{
\begin{array}{ll}
<0,\  &(L_1,L_2)\subset(-L^*,L^*),\\[2mm]
=0,\  &-L_1=L_2=L^*,\\[2mm]
>0,\  &(-L^*,L^*)\subset(L_1,L_2).
\end{array} \right.
\eess

\begin{thm}
\label{thm5.7}
Assume $0<a(0)<d$ holds, then

$(i)$ if $h_0\geq L^*$, then spreading always occurs for \eqref{a02}.

$(ii)$ if $h_0<L^*$, then there exists a positive constant $k^*$ such that $h_{\infty}-g_{\infty}=\infty$ when $k>k^*$.
\end{thm}

\bpf $(i)$ holds from Theorem \ref{thm5.2} since  $$
\lambda_{p}(\mathcal{L}_{\{(g_{\infty},\,h_{\infty}),\,d\}}+a(x))
>\lambda_{p}(\mathcal{L}_{\{(-L^{*},\,L^{*}),\,d\}}+a(x))=0.
$$

In what follows, we prove $(ii)$. Notice that
$$
-\mu_{2}+\beta(m(x),I,x)S-\gamma(b(x),I,x)>-\mu_{2}-\gamma(b(x),I,x)>-C
$$
for some $C>0$. Clearly $I(t,x)$ satisfies
\begin{equation*}
\left\{
\begin{array}{ll}
I_{t}\geq d\int_{g(t)}^{h(t)}J(x-y)I(t,y)dy-dI(t,x)-C I(t,x), &t>0,\, x\in(g(t), h(t)),\\[2mm]
I(t,x)=0,  &t\geq0,\ x\in\mathbb{R}\backslash(g(t),h(t)),\\[2mm]
h'(t)=k\int_{g(t)}^{h(t)}\int_{h(t)}^{+\infty}J(x-y)I(t,x)dydx, &t>0,\\[2mm]
g'(t)=-k\int_{g(t)}^{h(t)}\int_{-\infty}^{g(t)}J(x-y)I(t,x)dydx, &t>0,\\[2mm]
g(0)=-h_{0},\ h(0)=h_{0},   &x\in\mathbb{R},\\[2mm]
I(0,x)= I_0(x),   &x\in(-h_0,h_0),
\end{array} \right.
\end{equation*}
thereby, for any given constant $M$, there exists $k^*>0$ such that $h_{\infty}-g_{\infty}>M$
provided that $k>k^*$ by Lemma 3.9 in \cite{DWZ}. Then, $h_{\infty}-g_{\infty}=+\infty$ by
the arbitrariness of $M$.
\bigskip
\epf

Noting that the comparison principle for problem \eqref{a02} is not valid, we
cannot obtain the monotonicity of the solution for \eqref{a02} with $k$ and thus
cannot take $k$ as a sharp criterion for the spreading-vanishing dichotomy as in \cite{CDLL}.
However, recalling Remark \ref{rmk1.1} and $(ii)$ in Theorem \ref{thm5.7}, we have the
following result:

\begin{thm}
\label{thm5.9}
Suppose that $0<a(0)<d$ and $h_0<L^*$. For problem \eqref{a02}, there exists $0<k_{*}\leq k^{*}$ such that vanishing occurs if $0<k<k_{*}$ and spreading happens provided that $k>k^{*}$.
\end{thm}

Finally, we will discuss the impact of the diffusion coefficient on the vanishing and spreading of infectious disease.

Assume that $\int_{-\infty}^{-2h_0}J(z)dz>0$ (or $\int_{2h_0}^{+\infty}J(z)dz>0$) holds. Using Theorem \ref{thm3.3} with $L_1=-h_0$ and $L_2=h_0$, if $\max\limits_{x\in[-h_0,\,h_0]}a(x)<0$,
for any $d>0$, we have $\lambda_{p}(\mathcal{L}_{\{(-h_0,\,h_0),\,d\}}+a(x))<0$, while if $\max\limits_{x\in[-h_0,\,h_0]}a(x)>0$, there exists a $d^*>0$ such that
\bess
\lambda_{p}(\mathcal{L}_{\{(-h_0,\,h_0),\,d\}}+a(x))\left\{
\begin{array}{ll}
<0,\  &d>d^*,\\[2mm]
=0,\  &d=d^*,\\[2mm]
>0,\  &d<d^*.
\end{array} \right.
\eess

Inspired by the above analysis, combined with Theorem \ref{thm5.2}, we can obtain the following result.
\begin{thm}
\label{thm5.8}
Suppose $\int_{-\infty}^{-2h_0}J(z)dz>0$ (or $\int_{2h_0}^{+\infty}J(z)dz>0$) holds. Then,

$(i)$ if $\max\limits_{x\in[-h_0,\,h_0]}a(x)>0$, there exists a $d^*>0$ such that for $d<d^*$, then spreading occurs; conversely, if $d>d^*$, then vanishing occurs provided that  $\|S_0(x)\|_{L^{\infty}(\mathbb{R})}+\|I_{0}(x)\|_{C([-h_0,\,h_0])}$ is small enough;

$(ii)$ if $\max\limits_{x\in[-h_0,\,h_0]}a(x)\leq0$, for any $d>0$, then vanishing always occurs
as long as $\|S_0(x)\|_{L^{\infty}(\mathbb{R})}+\|I_{0}(x)\|_{C([-h_0,\,h_0])}$ is adequately small.
\end{thm}

\section{Discussion}
In this paper, we study a free boundary problem \eqref{a02} with media coverage and
hospital bed numbers, which describes a nonlocal diffusive SIS epidemic model. The free boundary
describes the moving front of the infected individuals, and the nonlocal diffusion operator
characterizes the long-distance spatial movement of individuals.

For the SIS model with nonlocal diffusion and free boundaries \eqref{a02}, the existence
and uniqueness of the global solution are given by using two fixed point theorems (see Theorem \ref{thm2.2}).
Then, we define the principal eigenvalue of the integral operator,
and analyze the impacts of media coverage and hospital bed number (Theorem \ref{thm3.2-1}),
interval length (Theorem \ref{thm3.2}), and diffusion coefficient (Theorem \ref{thm3.3})
on the principal eigenvalue. In addition, sufficient conditions
for disease spreading and vanishing (see Theorems \ref{thm5.2}, \ref{thm5.3}
and \ref{thm5.4}) are given. Finally,
we discuss the impact of the principal eigenvalue
on the spreading or vanishing of infectious diseases.
If $a(0)\geq d$, then $\lambda_{p}(\mathcal{L}_{\{(L_1,\,L_2),\,d\}}+a(x))>0$ for any $L_1<L_2$, and the disease is always spreading (see Theorem \ref{thm5.6}). If $0<a(0)<d$, there exists an $L^*>0$,
then spreading always appears for $h_0\geq L^*$ (see Theorem \ref{thm5.7}); and
when $h_0<L^*$, the impact of expanding capability $k$ on the spreading or vanishing
of disease is discussed. That is, there exists $0<k_{*}\leq k^{*}$
such that vanishing occurs if $0<k<k_{*}$, and spreading happens provided that $k>k^{*}$
(see Theorem \ref{thm5.9}). If $\max_{x\in[-h_0,\,h_0]}a(x)>0$, there exists a $d^*>0$
such that for $d<d^*$, then the disease spreads; if $d>d^*$ and $\|S_0(x)\|_{L^{\infty}(\mathbb{R})}+\|I_{0}(x)\|_{C([-h_0,\,h_0])}$
is small enough,
then the disease vanishes; if $\max_{x\in[-h_0,\,h_0]}a(x)\leq0$, then $d^*=0$, that
is, then vanishing always appears provided that $\|S_0(x)\|_{L^{\infty}(\mathbb{R})}+\|I_{0}(x)\|_{C([-h_0,\,h_0])}$
is adequately small (see Theorem \ref{thm5.8}).

Finally, we may conclude that
the differences between nonlocal diffusion in \eqref{a02} and local diffusion
in \eqref{a} in our mathematical analysis are as follows:
first, the existence and uniqueness
of the global solution for \eqref{a} are obtained by straightening the boundary and the
first-order fixed point theorem. However, owing to lack of compactness, the existence and uniqueness of global
solutions for \eqref{a02} are given by using two fixed point theorems. Second, for
\eqref{a}, the corresponding principal eigenvalue always exists. However, for the nonlocal diffusion
problem, the principal eigenvalue may not exist. In this paper, for the nonlocal
diffusion model \eqref{a02}, we first define the generalized principal
eigenvalue $\lambda_{p}(\mathcal{L}_{\{(L_1,\,L_2),\,d\}}+a(x))$ of the integral operator,
and then show that the generalized principal eigenvalue is the principal eigenvalue under condition $(\mathbf{H})$. Third,
unlike local diffusion, whose principal eigenvalue is clear, the nonlocal operator leads to more possibilities because of the choice of the kernel function.

It is worth mentioning that model \eqref{a02} incorporates media coverage and
hospital bed numbers. Based on the monotonicity of the generalized principal eigenvalue on media
coverage and hospital bed numbers, we study the influence of the principal eigenvalue on infectious diseases, which implies that
large media coverage and hospital bed numbers are beneficial
to the prevention and control of disease.

\end{document}